\documentclass[12pt,twoside]{article} 
 \input mfpic \opengraphsfile{rysChinSim09}
\usepackage{amssymb}
\usepackage{graphicx}
\usepackage{amsmath}%
\usepackage{amsfonts}

\newenvironment{proof}[1][Proof]{\textbf{#1.} }{\ \rule{0.5em}{0.5em}}
\def\Mod{{\operatorname{Mod}}}
\def\Comod{{\textrm{-}\operatorname{Comod}}}
\def\comod{{\textrm{-}\operatorname{comod}}}
\def\soc{{\operatorname{soc}}}
\def\End{{\operatorname{End}}}
 \def\Ext{{\operatorname{Ext}}}

 \def\Hom{{\operatorname{Hom}}}
  \def\hom{{\operatorname{hom}}}
 \def\inj{{\operatorname{\textrm{-}inj}}}
 \def\proj{{\operatorname{\textrm{-}proj}}}

 \def\rep{{\operatorname{rep}}}
 \def\bdim{{{\textbf{dim}}}}
 \def\len{{{{\textbf{lgth}}}\,}}
 \def\SQR{{\tiny\mbox{$\square$}}}
  \def\Car{{\Large{\mbox{$\mathfrak{c}$}}}}
  
 \def\mapr#1#2{\smash{\mathop{\longrightarrow}\limits^{#1}_{#2}}}
\def\mapl#1#2{\smash{\mathop{\longleftarrow}\limits^{#1}_{#2}}}
\def\mapupp#1#2{\llap{$\vcenter{\hbox{$\scriptstyle#1$}}$}\Big\uparrow
\rlap{$\vcenter{\hbox{$\scriptstyle#2$}}$}}
\def\mapdownn#1#2{\llap{$\vcenter{\hbox{$\scriptstyle{#1}$}}$}\Big\downarrow
\rlap{$\vcenter{\hbox{$\scriptstyle#2$}}$}}
\def\RightLeftarr#1#2{\vcenter{\vbox{\offinterlineskip\hbox to
1.7cm{\hfil$\scriptsize{#1}$\hfil}
\vskip0.5pt
\hbox to 1.8cm{\rightarrowfill}
 \vskip-0.7pt
\hbox to 1.8cm{\leftarrowfill}
\vskip0.5pt
\hbox to 1.7cm{\hfil$\scriptsize{#2}$\hfil}}}}

\def\shp{\hbox{\lower3pt\hbox{\vbox{
\mfpic[10][10]{-0.5}{1}{-0.5}{0.5}
\hatchwd{1}
\thatch[30,45]\ellipse{(0,0.2),1.5,1.0}
\ellipse{(0,0.2),1.5,1.0}
\endmfpic}}}}

\def\shk{\hbox{\lower3pt\hbox{\vbox{
\mfpic[10][10]{2.2}{1}{0}{0.5}
\shade[15]\ellipse{(0,0.5),1.5,1.0}
\ellipse{(0,0.5),1.5,1.0}
\endmfpic}}}} 
\begin{document}

\title{Coxeter transformation and     inverses of Cartan matrices  
 for coalgebras\footnote{$\mbox{}$\hspace{1.6em}Key words: coalgebra, comodule, pseudocompact algebra, Cartan matrix, Coxeter transformation, 
 quiver}\vspace{-0.3cm}}
\author{{\small \bf William Chin (DePaul University)}\\ {\small \bf and} \\
 {\small \bf Daniel Simson\footnote{Supported by Polish Research Grant  1 P03A N N201/2692/35/2008-2011} (Nicolaus Copernicus University)}}
\date{}
\maketitle
\vspace{-1.1cm}
 
\begin{abstract}
\smallskip

Let $C$ be a coalgebra and let 
$\mathbb{Z}^{I_C}_\blacktriangleright, \mathbb{Z}^{I_C}_\blacktriangleleft\subseteq \mathbb{Z}^{I_C}$ be  the Grothendieck groups of the category $C^{op}\inj$ and $C\inj$ of the socle-finite injective right and left $C$-comodules, respectively. One  of the main aims of the paper is to study Coxeter transformation $\mathbf{\Phi}_C:\mathbb{Z}^{I_C}_\blacktriangleright\,\to\,  \mathbb{Z}^{I_C}_\blacktriangleleft$   and  its dual    $  \mathbf{\Phi}^-_C: \mathbb{Z}^{I_C}_\blacktriangleleft \to   \mathbb{Z}^{I_C}_\blacktriangleright$ of a pointed   sharp Euler coalgebra $C$, and to relate  the action of $ {\mathbf{\Phi}}_C$ and $ {\mathbf{\Phi}}^{-}_C$ on  a  class  of  indecomposable finitely cogenerated $C$-comodules $N$ with the ends of  almost split sequences starting with $N$  or  \hbox{ending at $N$}. By applying \cite{CKQ}, we also show that    if  $C$ is a  pointed  $K$-coalgebra such that the  every vertex of the left Gabriel quiver    ${}_CQ$  of $C$ has only finitely many neighbours   then
for any   indecomposable non-projective left $C$-comodule $N$ of finite $K$-dimension,   there exists a unique almost split sequence $  
0  \, \mapr{}{} \,\  \tau_C N   \,\mapr{}{}  \,N' \,\mapr{}{} \,N  \,
\mapr{}{} \, 0  $ 
in  the category $C\Comod_{fc}$ of finitely cogenerated left $C$-comodules,  with an \hbox{indecomposable comodule $\tau _CN$}.  We show that $\mathbf{dim}\, \mathbf{\tau}_C N =\mathbf{\Phi}_C(\mathbf{dim}\,N)$, if $C$ is hereditary, or more generally, if $\textrm{inj.dim} \, DN =1$ and $ \textrm{Hom}_C(C, DN)=0$. \end{abstract}

\markboth{\hfill\sc w. chin and d. simson\hfill}{\hfill\sc cartan matrix and   coxeter transforation for coalgebras
\hfill }
\bigskip

 Throughout we fix an arbitrary  field $K$ and   $D(-) = (-)^*= \textrm{Hom}_K(-, K)$ is the ordinary $K$-linear duality functor. We recall that  a $K$-coalgebra $C$ is said to be pointed  if all simple $C$-comodules are one-dimensional.  Let $C$ be a pointed   $K$-coalgebra and $C^* = \textrm{Hom}_K(-, K)$ the $K$-dual (pseudocompact \cite{Sim2},  \cite{Sim3}) $K$-algebra with respect to the convolution product, see \cite{DNR}, \cite{Montg}.   
 We denote by $C\mbox{-}\textrm{Comod}$ and $C\mbox{-}\textrm{comod}$ the category of left $C$-comodules and finite-dimensional left $C$-comodules, respectively. The corresponding categories of right $C$-comodules are denoted by   $C^{op}\mbox{-}\textrm{Comod}$ and $C^{op}\mbox{-}\textrm{comod}$. The socle of
a comodule $M$ in $ C\Comod$ is denoted by $ \soc\,M$

We recall from \cite{CKQ} and \cite {CKQ2} that, for a class of coalgebras $C$ (including left semiperfect ones), given an indecomposable non-injective $C$-comodule $M$ in $C\comod$ and an indecomposable non-projective $C$-comodule $N$ in $C\comod$ there exist  almost split sequences \vspace{-0.2cm}
$$  
0  \, \mapr{}{} \,\   M   \,\mapr{}{}  \,M' \,\mapr{}{} \,\tau_C^{-}M  \,
\mapr{}{} \, 0\quad\mbox{and}\quad
0  \, \mapr{}{} \,\  \tau_C N   \,\mapr{}{}  \,N' \,\mapr{}{} \,N  \,
\mapr{}{} \, 0   \vspace{-0.2cm} \leqno(*)
$$ 
in $C\Comod$, where 
     $
  C\mbox{-} {\textrm{Comod}}^\bullet_{fc}   \, \,\RightLeftarr{\tau^-_C}{\tau_C}\, \,  C\mbox{-}{\textrm{comod}}_{f\mathcal{P}}  $ are the Auslander-Reiten translate  operators   (1.17).  
 On the other hand, for a class of computable coalgebras $C$, a Cartan matrix $\Car_C\in \mathbb{M}_{I_C}(\mathbb{Z})
$, its inverse ${\Large{\mbox{$\mathfrak{c}$}}}_C^{-1}$,   and  a corresponding Coxeter transformation \hbox{$ \mathbf{\Phi}_C\!:\! K_0(C)\! \to\! K_0(C) $}  is defined and studied in  \cite{Sim4} (see also \cite{Ch} and \cite{Gr}), where $K_0(C)= K_0(C\comod)\cong \mathbb{Z}^{(I_C)}$ is the Grothendieck group of $C\comod$.  

One of the main aims of this paper is to construct an inverse $\Car_C^{-1}$ (left and right) of the Cartan matrix $\Car_C$ and  Coxeter transformations $\mathbf{\Phi}_C:\mathbb{Z}^{I_C}_\blacktriangleright\,\to\,  \mathbb{Z}^{I_C}_\blacktriangleleft$, $  \mathbf{\Phi}^-_C: \mathbb{Z}^{I_C}_\blacktriangleleft \to   \mathbb{Z}^{I_C}_\blacktriangleright$, for any computable    coalgebra $C$ such that any simple left (and right) $C$-comodule admits a finite and socle-finite injective resolution. We prove that, under a suitable assumption on indecomposable $C$-comodules  $N$ and $M$, there exist  almost split sequences $(*)$  in $C\Comod$ and the following equalities hold (compare with \cite[Corollary IV.2.9]{ASS} and \cite[pp.67]{Sim4})\vspace{-0.2cm}
$$
\mathbf{dim}\, \mathbf{\tau}_C^{-}(M)=\mathbf{\Phi}_C^{- }(\mathbf{dim}\, M)\quad\textrm{and} \quad   \mathbf{dim}\, \mathbf{\tau}_C(\,N)=\mathbf{\Phi}_C(\mathbf{dim}\,N),\vspace{-0.2cm}
$$
where $\mathbf{dim}\, X \in  \mathbb{Z}^{I_C}$ is   the dimension vector  of the comodule $X$, see Section 2.  

We recall from \cite{KR} that a coalgebra  $C$ is said to be  left locally artinian if every indecomposable injective left $C$-comodule is artinian. Recall  also that  a  coalgebra over an algebraically closed   field   is  pointed if and only if it is basic, see \hbox{\cite{Ch}, \cite{CMo}, \cite[p.404]{Sim0}, \cite[5.5]{Sim1},  \cite[2.2]{Sim5}}.  
  \vspace{-0.4cm}
 
\begin{center}
\section{Preliminaries on comodule categories}
\end{center}
\vspace{-0.2cm}

Let $C$ be a $K$-coalgebra. We collect in this section basic  facts concerning $C$-comodules, pseudocompact  $C^*$-modules,  the existence of almost split sequences in $C\Comod$, 
duality and injectives in the category of comodules. 

We recall from \cite{Ch}, \cite{DNR}, \cite{Montg} and \cite{Sim1} that any left $C$-comodule $M$ is viewed as a rational (=discrete) right module over the pseudocompact algebra $C^*$ and $M^*=  D(M) = \Hom_K(M,K)$ is a  pseudocompact  left $C^*$-module. The functor $D(-)$ defines a duality\linebreak  $\widetilde D: C\Comod\,\,\mapr{}{}\,\, C^*\textrm{-PC}$, where $C^*\textrm{-PC}$ is the category of pseudocompact left $C^*$-modules. The quasi-inverse is the functor $(-)^\circ = \textrm{hom}_K(-, K)$ that associates to any  $Y$ in $ C^*\textrm{-PC}$ the left $C$ comodule $Y^\circ = \textrm{hom}_K(Y, K)$ consisting of all    continuous $K$-linear maps $Y \to K$.\linebreak 
 It follows from \cite{Ga73} that  the algebra $C^*$ is left (and right)  topologically semiperfect, that is, every simple left $C^*$-module admits a projective cover in $ C^*\textrm{-PC}$ (see also \cite{Sim1}); equivalently,  $ C^*$ admits a decomposition  
 $ C^*   \cong \prod \limits  _{j\in I } \Lambda  e_j   $
in $ C^*\textrm{-PC}$, where  $\{e_j\}_{j\in I}$ is a topologically complete  set of  pairwise orthogonal  primitive idempotents   such that $e_j\Lambda  e_j$ is a local algebra, for every $i\in I$. The decomposition  is unique up to isomorphism and  permutation.  
 
The coalgebra  $C$ 
(or more generally,    any $C$-$C$-bicomodule)  can be viewed as a bimodule over the algebra $C^*$ with respect to the  right and the left hit actions of   $C^*$ on $C$, usually denoted by the
symbols $\leftharpoonup,\rightharpoonup$ as in  \cite{DNR} and  \cite{Montg}. Here we  omit these symbols and simply use juxtapostion, e.g., $eC
=e\rightharpoonup C$ and $ Ce
=  C\leftharpoonup e$, for any $e\in C^*$. Notice that $C e$ is an injective
right $C$-comodule and $eC $ is an injective
left $C$-comodule, for any idempotent $e\in C^*$.
  
  The following two simple lemmata are often used in the paper.  \medskip

 {\sc Lemma 1.1.}  \textit{Assume that $C$ is a coalgebra, $e=e^{2}$ is an idempotent in $ C^*$, and $D(-)= \operatorname*{Hom}_K(-,K)$.}

 (a)  \textit{There is an isomorphism $\widetilde D( eC)\cong  C^*e$ of  left
$C^*$-modules.}

 (b)  \textit{If $C$ 
is of finite dimension, then there is an isomorphism
$  D(C^*e)\cong eC$ of  right $C^*$-modules.}

 (c)  \textit{$\operatorname*{Hom} 
_{C}(C',C e)=\operatorname*{Hom}_{C}(C',C' e)$
for every subcoalgebra $C'$ of $C$.}
\smallskip

\begin{proof}  See \cite{CKQ2} and  \cite{CG}.
\end{proof}\smallskip

 {\sc Lemma 1.2.}  \textit{Let $C$ be a $K$-coalgebra. Given a left comodule $M$ in} $C\comod$,   \textit{the $K$-dual space} $D(M) = \Hom_K(M,K)$  \textit{admits a natural structure of right $C$-comodule and}  $D(-) = \Hom_K(-,K)$ {defines the pair of  dualities }\vspace{-0.2cm}
 $$
  C\comod   \,\, \RightLeftarr{D}{D}\,\, C^{op}\comod\vspace{-0.2cm} \leqno(1.3)
 $$  
 
\begin{proof}  See \cite{Ch}, \cite{Gr},  \cite{Lin},   \cite{Sim1},  \cite{Tak}.
\end{proof}\medskip

An \textbf{injective copresentation} of a comodule $M$ in $C\Comod$ is an exact sequence \vspace{-0.2cm}
$$
0\, \mapr{}{}\,M\, \mapr{}{}\,E_{0}\, \mapr{g}{}\, E_{1} , \vspace{-0.2cm}\leqno(1.4)
$$
where $E_{0}$ and $E_{1}$ are
injective comodules. We call a comodule $M$ in $C\Comod$ (socle) \textbf{finitely copresented}  if $M$ admits a socle-finite injective copresentation, that is, the injective comodules $E_0$ and $E_1$ have finite-dimensional socle.  We denote by $C\mbox{-}\textrm{Comod}_{fc}$ the full subcategory of $C\mbox{-}\textrm{Comod}$ whose objects are the finitely copresented comodules, and  by $C\mbox{-}\textrm{inj}$  the full subcategory of $C\mbox{-}\textrm{Comod}_{fc}$ whose objects are the socle-finite   injective  comodules. We set $C\comod_{fc}= C\comod\cap C\Comod _{fc}$. Finally, we denote by $C\mbox{-}\overline{\textrm{Comod}}_{fc}= \textrm{Comod}_{fc}/\mathcal{I}$ the quotient category of $C\mbox{-}\textrm{Comod}_{fc}$ modulo the two-sided ideal $\mathcal{I}=[C\mbox{-}\textrm{inj}]$     consisting of all $f\in \Hom_C(N, N')$, with $N$ and $N'$ in $C\Comod_{fc}$,  that have a factorisation through a socle-finite injective comodule. 

It is observed in \cite{KR} that $C\Comod_{fc}$ is an abelian category if and only if $C$ is left cocoherent. In this case $C\Comod_{fc}$ is closed under extensions in $C\Comod$ and  contains minimal injective resolutions of comodules $M$ in $C\Comod_{fc}$, see \cite[Section 3]{Sim4}.

We recall that a comodule $M$ is \textbf{quasi-finite} if $\dim_K  \textrm{Hom}_C(X, M)$ is finite, for any $X$ in $C\comod$;   equivalently, if
the simple summands of $\soc\, M$ have finite (but perhaps unbounded)
multiplicities \cite{Ch}, \cite{Tak}. It is easy to check that every socle-finite comodule is quasi-finite. Hence all comodules in $C\mbox{-}\textrm{Comod}_{fc}$ are quasi-finite.

  Given a left quasi-finite $C$-comodule $M$, the covariant  \textbf{cohom functor} \vspace{-0.2cm}
  $$
  h_C(M, -):C\Comod \,\,\mapr{}{}\,\, \Mod(K)\vspace{-0.2cm}
  $$ 
   is defined by associating to any comodule $N$ in $C\Comod$ the vector space  $h_C(M, N) = {\lim\limits_{\longrightarrow _\lambda}} D\Hom_C(N_\lambda, M)$, where $\{N_\lambda\}$ is the family of all finite-dimensional subcomodules \hbox{of $M$ \cite{Tak}}. 

 Denote by $C^{op}\Comod_{fp}$   the full subcategory of $C^{op}\Comod$ whose objects are the (injectively)  finitely presented $C^{op}$-comodules, that is, the $C^{op}$-comodules  $L$ that admit a short exact sequence $E'_1\, \mapr{g'}{}\, E'_0 \, \mapr{}{}\, L\, \mapr{}{}\, 0$ in $C^{op}\Comod$, with socle-finite injective comodules $E'_1$ and $E'_0$, called a socle-finite injective presentation of $L$. 
  Following \cite[Section 3]{CKQ}, we define a pair of  contravariant left exact functors   \vspace{-0.2cm}
 $$
  C\Comod_{fc} \,\,\RightLeftarr{ \nabla_C}{ \nabla'_C}\,\, C^{op}\Comod_{fp}\vspace{-0.2cm}\leqno(1.5)
 $$
  to be the composite functors making the following diagrams commutative\vspace{-0.2cm}
  $$
  \begin{array}{ccc}
 C\Comod_{fc} &\,\, \mapr{\widetilde D}{\simeq}\,\,&C^*\textrm{-PC}_{fp}\vspace{1ex}\\
  \mapdownn{ \nabla_C}{}
  && \mapdownn{(-)^+}{}\vspace{1ex}\\
 C^{op}\Comod_{fp}& \,\, \mapl{(-)^{\circ}}{\simeq}\,\, &{C^*}^{op}\textrm{-PC}_{fc},
  \end{array} \qquad
    \begin{array}{ccc}
 C\Comod_{fc} &\,\, \mapl{(-)^{\circ}}{\simeq} \,\,&C^*\textrm{-PC}_{fp}\vspace{1ex}\\
  \mapupp{\nabla'_C}{}
  && \mapupp{(-)^+}{}\vspace{1ex}\\
 C^{op}\Comod_{fp}& \,\, \mapr{\widetilde D}{\simeq}\,\, &{C^*}^{op}\textrm{-PC}_{fc},
  \end{array}
  \vspace{-0.2cm}\leqno(1.6)
 $$ 
  where  $C^*\textrm{-PC}_{fp}$ (resp. $ {C^*}^{op}\textrm{-PC}_{fc}$)  is the category of pseudocompact   (top-) finitely presented (resp.  (top-)  finitely copresented) modules (see \cite{DNR}, \cite{Sim2}, \cite{Sim3}),  $\widetilde D = \Hom_K(- , K)$,  \vspace{-0.2cm}
  $$
  (-)^+=  \textrm{hom}_{C^*}(-, C^*):C^*\textrm{-PC}_{fp}  \,\, \mapr{ }{}\,\,{C^*}^{op}\textrm{-PC}_{fc}\vspace{-0.2cm}
  $$
    is a contravariant  functor that associates to any $X$ in $C^*\textrm{-PC}_{fp}$,  with the top-finite pseudocompact projective presentation $P_1\, \mapr{f_1}{}\, P_0 \, \mapr{}{}\, X \, \mapr{}{}\, 0$, where $P_1, P_0$ are  finite direct sums of indecomposable projective $C^*$-modules, the right $C^*$-module $X^+ = \textrm{hom}_{C^*}(X, C^*)$ of all continuous $C^*$-homomorphisms $X \to C^*$, with the   top-finite pseudocompact projective copresentation \vspace{-0.2cm}
    $$
    0  \, \mapr{}{}\,  X^+ \, \mapr{}{}\, P^+_0\, \mapr{f^+_1}{}\, P_1^+ .  \vspace{-0.2cm}
    $$
     Finally, $Y^\circ = \textrm{hom}_K(Y, K)$ consists of all continuous $K$-linear maps $Y \to K$ and $(-)^\circ$ associates to $X^+$ the right $C$-comodule $(X^+)^\circ$ in $ C^{op}\Comod_{fp}$, with the socle-finite injective presentation 
     \vspace{-0.2cm}
     $$
     (P_1^+)^\circ\, \mapr{(f_1^+)^\circ}{}\, ( P_0^+)^\circ \, \mapr{}{}\, (X^+)^\circ \, \mapr{}{}\, 0. \vspace{-0.2cm}
     $$
     The functors in the right hand diagram of (1.6) are defined analogously. Sometimes, for simplicity of the notation, we write $\nabla_C$ instead of $\nabla'_C$.
  
  Following  \cite{CKQ} and the classical construction of Auslander \cite{Ausl}, we define  the  Auslander  \textbf{transpose operator} \vspace{-0.25cm}
  $$
  \textrm{Tr}=  \textrm{Tr}_C:
  C \textrm{-Comod}_{fc}  \, \mapr{}{}\,C^{op} \textrm{-Comod}_{fc}\vspace{-0.25cm}\leqno(1.7)
  $$
   (on objects only!) that associates to any comodule  $M$ in $C \textrm{-Comod}_{fc}$, with a minimal socle-finite injective copresentation    (1.4),   the comodule \vspace{-0.25cm}
   $$
   \textrm{Tr}_C M = \textrm{Ker} [ \nabla_C E_1 \, \mapr{ \nabla_C(g)}{}\,  \nabla_C E_0 ] \vspace{-0.25cm}
   $$
    in $C^{op} \textrm{-Comod}_{fc} $.  Basic properties of $\textrm{Tr}_C$ are listed in  \cite[Proposition 3.2]{CKQ}.
    
  The existence of almost split sequences in $C\Comod_{fc}$ essentially depends on the following theorem  slightly extending  some of the results in  \cite{CKQ} and  \cite{CKQ2}.\medskip
  
 {\sc Theorem 1.8.}  \textit{Let $C$ be a $K$-coalgebra and $ \nabla_C$ the functor} \textrm{(1.5)}.

\textrm{(a)}  \textit{There are functorial isomorphisms}  $ \nabla_C M  \cong \Hom_C( C, M)^\circ  \cong  h_C( M, C)$, \textit{for any comodule $M$ in} $C \textrm{-Comod}_{fc}$.

 \textrm{(b)}  \textit{The functors $ \nabla_C$, $ \nabla'_C$ are left exact and   restrict  to the  dualities} \vspace{-0.3cm}
 $$
  C\mbox{-}\textrm{inj}\,\, \RightLeftarr {\nabla_C}{\nabla_C}\,\, C^{op}\mbox{-} \textrm{inj}\vspace{-0.3cm}
  \leqno(1.9)
  $$
    that are quasi-inverse to each other.  \textit{Moreover, given an idempotent $e\in C^*$, the comodule $Ce$ lies in} $ C^{op}\mbox{-}\textrm{inj}$,   \textit{the comodule} $eC$  \textit{lies in} $ C\mbox{-}\textrm{inj}$,  \textit{and   there is an isomorphism $  \nabla_C(C e) \cong eC$ of left $C$-comodules}.

 \textrm{(c)}  \textit{For any  comodule $M$ in} $C \textrm{-Comod}_{fc} $,  \textit{with a minimal socle-finite injective copresentation}    (1.4),  \textit{the comodules} $\textrm{Tr}_CM $, $ \nabla_C E_1 $, $  \nabla_C E_0 $  \textit{lie in}  $C^{op} \textrm{-Comod}_{fc}$, $  \nabla_C M $  \textit{lies in}  $C^{op} \textrm{-Comod}_{fp}$,  \textit{and the following sequence}\vspace{-0.2cm}
   $$
  0  \, \mapr{}{}\,   \textrm{Tr}_C M    \, \mapr{}{}\,  \nabla_C E_1 \, \, \mapr{ \nabla_C  (g)}{}\,  \, \nabla_C E_0   \, \mapr{}{}\,  \nabla_C M  \, \mapr{}{}\,0\vspace{-0.2cm} \leqno(1.10)
  $$
 \textit{ is exact in} $C^{op} \textrm{-Comod}$.

  \textrm{(d)}  \textit{The transpose operator } $\textrm{Tr}_C $, {\it together with the functor $ \nabla_C$, induces the equivalence of quotient categories}  $ \textrm{Tr}_C  :  C\mbox{-}\overline{\textrm{Comod}}_{fc}   \, \, \mapr{\simeq}{}\, \,  C^{op}\mbox{-}\overline{\textrm{Comod}}_{fc}$.  \medskip
 
\begin{proof}  For our future purpose and the convenience of the reader,  we   outline   the proof. 

(a) Let $\{C_\lambda\}$ is the family of all finite-dimensional subcoalgebras of $C$ and let $M$ be a comodule in $C \textrm{-Comod}_{fc}$. Then $M$ is quasi-finite, $C  \cong {\lim\limits_{\longrightarrow _\lambda}}C_\lambda$,   and we get isomorphisms  \vspace{-0.3cm}
$$ 
 \begin{array}{rcl}
 \nabla_CM= ((\widetilde DM)^+) ^\circ & \cong &\hom_{C^*}(\widetilde DM, C^*) ^\circ 
\\
& \cong &
   \Hom_{C}(C, M) ^\circ  \\
& \cong &
[\lim\limits_{\longleftarrow _\lambda} \Hom_{C}(C_\lambda, M) ]^\circ   
\\
& \cong & \lim\limits_{\longrightarrow _\lambda}  \Hom_{C}(C_\lambda, M) ^\circ 
\\
& \cong & \lim\limits_{\longrightarrow _\lambda}  D\Hom_{C}(C_\lambda, M) 
\\
& =&  h_C(M, C).
\end{array}\vspace{-0.2cm}
 $$
 One can easily see that the composite isomorphism is functorial at $M$.

(b) Apply the definition of $ \nabla_C$.

To prove  (c)  and 
(d), we note that  the exact functors $\widetilde D:
 C\Comod_{fp} \,\, \mapr{ }{}\,\,C^*\textrm{-PC}_{fp}$ and $(-)^\circ :{C^*}^{op}\textrm{-PC}_{fc} \,\, \mapr{}{}\,\,  C^{op}\Comod_{fc}  $ defining the functor $ \nabla_C$ are equivalences of categories carrying injectives to projectives and projectives to injectives, respectively.  Recall    that $C^*$ is a topological semiperfect algebra. Now, given an indecomposable   comodule $M$ in $C \textrm{-Comod}_{fc} $, with a minimal socle-finite  injective copresentation (1.4), we get a    pseudocompact minimal top-finite projective presentation \vspace{-0.2cm}
 $$
 \widetilde DE_1\, \mapr{\widetilde Dg}{}\, \widetilde DE _0 \, \mapr{}{}\,\widetilde DM  \, \mapr{}{}\, 0,\vspace{-0.2cm}
 $$
in $ C^*\mbox{-}\text{PC}$,  with $\widetilde DE_1= E^*_1$, $\widetilde DE_0= E^*_0$ finite direct sums of indecomposable projective $C^*$-modules, of  the right pseudocompact $C^*$-module $\widetilde DM$. Hence, by applying the left exact functor  $\textrm{hom}_{C^*}(-  , C^*)$,   and the definition of the Auslander transpose $\textrm{Tr}_{C^*}(\widetilde DM)$ of the pseudocompact left $C^*$-module $\widetilde DM$, we get    the exact sequence \vspace{-0.2cm}
    $$
  0  \, \mapr{}{}\,   ( \widetilde DM )^+   \mapr{}{}\,  (\widetilde DE_0)^+ \, \mapr{(\widetilde Dg )^+ }{}\,  (\widetilde DE_1) ^+\, \mapr{}{}\,  \textrm{Tr}_{C^*}(\widetilde DM)  \, \mapr{ }{}\, 0\vspace{-0.2cm}\leqno(1.11)
    $$
in ${C^*}^{op}\textrm{-PC}$ and the  projective copresentation
 $0  \, \to\,  ( \widetilde DM )^+ \mapr{}{}\,  (DE_0)^+  \mapr{(\widetilde Dg )^+ }{}  (DE_1) ^+$ of  the  right  pseudocompact $C^*$-module $(\widetilde DM )^+$, where  $ (\widetilde D E_0)^+$ and $ (\widetilde DE_1)^+$ are finitely generated projective top-finite right $C^*$-modules.  The sequence (1.11) induces the sequence (1.10) and (c) follows.
 The statement (d) follows from the corresponding properties of the Auslander transpose operator $\textrm{Tr}_{C^*}: C^*\textrm{-PC}_{fp}\,\mapr{}{}\, {C^*}^{op}\textrm{-PC}_{fp}$    on  the pseudocompact  
 finitely presented top-finite modules over $C^*$, see \cite[Proposition IV.2.2]{ASS}, \cite[Section IV.1]{ARS}, \cite[Proposition 11.22]{Si92}.
\end{proof}\medskip

We denote by  $C\Comod^\bullet_{fc}$ and by $C\Comod^\nu_{fc}$  the full subcategory of $C\Comod_{fc}$ consisting of the comodules $M$ such that  $\dim_K\textrm{Tr}_C(M)$ is finite and $\dim_K (\widetilde DM)^+$ is finite, respectively.  Following the representation theory of finite-dimensional algebras, we define the 
\textbf{Nakayama functor}  (covariant)  \vspace{-0.2cm}
$$
\nu_C : C\Comod^\nu_{fc} \, \mapr{}{}\,  C\comod\vspace{-0.2cm}\leqno(1.12)
$$
 by the formula $\nu_C ( -)=D \nabla_C(-)$.

A coalgebra C is said to be \textbf{left semiperfect} \cite{Lin} if every simple left
comodule has a projective cover, or equivalently, the injective envelope $E(X)$ of any finite-dimensional right  $C$-comodule $X$ is   finite-dimensional.

It is easy to see that, for a left semiperfect coalgebra $C$, the functor $\nu_C$ restricts to the equivalence of categories\vspace{-0.3cm}
$$
\nu_C : C\inj \, \mapr{\simeq}{}\,  C\proj,\vspace{-0.2cm}\leqno(1.13)
$$
where $C\proj$ is the category of top-finite projective comodules in $C\comod$.  

We denote by $C\comod_{f\mathcal{P}}$ the full subcategory of $C\comod$ consisting of the left  comodules $N$  that, viewed as   rational right ${C^*}$-modules,   have a minimal  top-finite projective presentation  $P_1 \to P_0 \to N \to 0$ in  ${C^*}^{op}\textrm{-PC}= \textrm{PC-}C^*$, that is,  $P_0$ and $P_1$ are top-finite projective modules in $ \textrm{PC}\mbox{-}C^* $. Here we make the identification \vspace{-0.2cm}
$$
C\mbox{-} {\textrm{comod}} \equiv  {\textrm{rat}\mbox{-}}C^* = {\textrm{dis}\mbox{-}}C^* \subseteq  \textrm{PC}\mbox{-}C^*  ,  \vspace{-0.2cm}
$$
in the notation of \cite[Section 4]{Sim1}, where $ {\textrm{rat}\mbox{-}}C^*$ is the category of finite dimenional rational right ${C^*}$-modules.

Finally, we denote by  \vspace{-0.2cm}
$$
C\mbox{-}\underline{\textrm{comod}}_{f\mathcal{P}}= C\mbox{-}{\textrm{comod}}_{f\mathcal{P}}/\mathcal{P} \vspace{-0.2cm}
$$
 the quotient category of $C\mbox{-}\textrm{comod}_{f\mathcal{P}}$ modulo the two-sided ideal $\mathcal{P}$ of $C\mbox{-}{\textrm{comod}}_{f\mathcal{P}}$  consisting of all $f\in \Hom_C(N, N')$, with $N$ and $N'$ in $C\comod_{f\mathcal{P}}$,  that have a factorisation through a projective right ${C^*}$-module, when  $f$ is viewed as a ${C^*}$-homomorphism  between the  rational right ${C^*}$-modules $N$ and $N'$.  

 If $C$ is left semiperfect then, in view of the exact sequence (1.10) in  $ C^{op}\Comod$, we have  $C\comod_{f\mathcal{P}}=  C\comod$, $C\Comod^\bullet_{fc}=  C\Comod_{fc}$,   $C\Comod^\nu_{fc}=  C\Comod_{fc}$ and, by applying $\nu_C$   to the sequence (1.10) we get the exact sequence 
  \vspace{-0.3cm}
    $$
  0  \, \mapr{}{}\,   \nu_C(M)   \mapr{}{}\,   \nu_C(E_0)  \mapr{\nu_C(g )  }{}\,  \nu_C(E_1)\, \mapr{}{}\,  D\textrm{Tr}_{C}(M)  \, \mapr{ }{}\, 0\vspace{-0.2cm}\leqno(1.14)
    $$
in    $ C\comod$.

 The following simple lemma is of importance.
 \medskip
 
   {\sc Lemma 1.15.} \textit{ Let $C$ be a pointed   $K$-coalgebra and let   ${}_CQ$ be the left Gabriel quiver  of $C$}.  
    
  \textrm{(a)}   \textit{The duality} $D   :  C\mbox{-} {\textrm{comod}}   \, \, \mapr{}{}\, \,  C^{op}\mbox{-}{\textrm{comod}}$   \textrm{(1.3)}
  \textit{restricts to the duality}\vspace{-0.3cm}
  $$
   D   :  C\mbox{-} {\textrm{comod}}_{f\mathcal{P}}   \, \, \mapr{}{}\, \,  C^{op}\mbox{-}{\textrm{comod}}_{fc } = C^{op}\mbox{-}{\textrm{comod}}\cap C^{op}\mbox{-}{\textrm{Comod}}_{fc }  ,\vspace{-0.3cm}
      $$
\textit{In particular, a left $C$-comodule $N$ lies in $C\mbox{-} {\textrm{comod}}_{f\mathcal{P}} $ if and only if the right $C$-comodule $D(N)$ is finitely copresented}.

 \textrm{(b)}  \textit{The following four conditions are equivalent:}
 
 \textrm{(b1)}  \textit{the   equality} $C\mbox{-} {\textrm{comod}}_{f\mathcal{P}} = C\mbox{-} {\textrm{comod}}$ \textit{holds,}
 
  \textrm{(b2)}    \textit{the   inclusion } $C^{op}\mbox{-} {\textrm{comod}}\subseteq  C^{op}\mbox{-} {\textrm{Comod}_{fc}}$ \textit{holds,}

    \textrm{(b3)} \textit{every simple  comodule in $C^{op}\mbox{-} {\textrm{Comod}}$ is finitely copresented},

  \textrm{(b4)}  \textit{the quiver}  ${}_CQ$  \textit{is right locally bounded, that is,  for every vertex $a$ of}   ${}_CQ$  \textit{there is only a finite number of arrows $a\,\mapr{}{}\, j$ in} ${}_CQ$.

 \textrm{(c)}  \textit{If $C$ is right locally artinian then  the   equality} $C\mbox{-} {\textrm{comod}}_{f\mathcal{P}} = C\mbox{-} {\textrm{comod}}$ \textit{holds.}
 
  \medskip
 
 \begin{proof} (a) Since we make the identification 
   $
C\mbox{-} {\textrm{comod}} \equiv  {\textrm{rat}\mbox{-}}C^* = {\textrm{dis}\mbox{-}}C^* \subseteq  \textrm{PC}\mbox{-}C^*   $
(in the notation of \cite[Section 4]{Sim1}),   there is a commutative diagram 
 $$
  \begin{array}{ccccc}
 C\comod &\,\, \mapr{id}{ }\,\,&{\textrm{dis}\mbox{-}}C^*&\subseteq &\textrm{PC}\mbox{-}C^* \vspace{1ex}\\
  \mapdownn{ D}{\cong}
  && \mapdownn{ (-)^\circ}{\cong}
  &&  \mapdownn{(-)^\circ}{\cong}\vspace{1ex} \\
 C^{op}\comod & \,\, \mapr{id}{}\,\, &C^*{\mbox{-}}\textrm{dis} &\subseteq &C^{op}\Comod,
  \end{array}   \vspace{-0.2cm} 
 $$  
 Then  a left $C$-comodule $N$ lies in $C\mbox{-} {\textrm{comod}}_{f\mathcal{P}} $ if and only if  there is an exact sequence $P_1 \to P_0 \to N \to 0$ in $ \textrm{PC}\mbox{-}C^* $, where $P_0$ and $P_1$ are top-finite projective modules in $ \textrm{PC}\mbox{-}C^* $, or equivalently,    $N$ lies in ${C^*}^{op}\textrm{-PC}_{fp}= \textrm{PC}_{fp}\mbox{-}C^*$. By applying the duality $(-)^\circ :{C^*}^{op}\textrm{-PC} \to C^{op}\mbox{-} {\textrm{Comod}}$, see (1.6),  we get an exact sequence $0\to N^\circ \to  P^\circ_0 \to P^\circ_1  $ in $C^{op}\mbox{-} {\textrm{Comod}}$. Since $\dim_K N$ is finite, we have  $ N^\circ=  D(N)$. This shows that  $D(N)$ lies in $C^{op}\mbox{-}{\textrm{comod}}_{fc } $, because $ P^\circ_0$ and $P^\circ_1  $ are socle-finite injective right $C$-comodules.  It follows that the   duality  (1.3) restricts to the duality
 $
   D   :  C\mbox{-} {\textrm{comod}}_{f\mathcal{P}}   \, \, \mapr{}{}\, \,  C^{op}\mbox{-}{\textrm{comod}}_{fc } $.
   
   (b) By (a),  the equality  $C\mbox{-} {\textrm{comod}}_{f\mathcal{P}} = C\mbox{-} {\textrm{comod}}$  holds if and only if the equality $C^{op}\mbox{-}{\textrm{comod}}_{fc } = C^{op}\mbox{-}{\textrm{comod}}$ holds, that is,    the conditions  (b1) and (b2)  are equivalent. 
   
 The  implication (b2)$\Rightarrow$(b3) is obvious.  To prove the inverse implication (b3)$\Rightarrow$(b2), we  assume that  the simple right $C$-comodules lie  in $C^{op} \textrm{-Comod} _{fc}$ and let $X$ be a comodule in $C^{op}\comod$. By standard arguments and the  induction on the $K$-dimension of $X$, we show that $X$ lies in $C\Comod_{fc}$ (apply the diagram in \cite[p. 13]{CKQ}).

   (b3)$\Rightarrow$(b4)
   Fix a direct sum decomposition $\soc \, C_C = \bigoplus _{j\in I_C} \widehat{S}(j)$ of the right socle $\soc\, C_C $ of $C$, where $I_C$ is an  index set and $\{\widehat S(j)\}_{j\in I_C}$ is a set of pairwise non-isomorphic simple right $C$-coideals. Denote by  $\widehat E(j) =E(\widehat S(j))$ the injective envelope of $\widehat S(j)$.
 
 It follows from the dualities (1.9) and \cite[Theorem 2.3(a)]{Sim06} that  the quiver ${}_CQ$ is dual to the  right Gabriel quiver ${}Q_C$ of $C$.   Hence, by the assumption (b2), for every vertex $a$ of the  quiver $ Q_C$, there is only a finite number of arrows $j\,\mapr{}{}\, a$ in $Q_C$. In other words, $\dim_K \Ext_C^1(\widehat S(j), \widehat S(a))$ is finite, for all $j\in I_C$, and $\dim_K \Ext_C^1(\widehat S(j), \widehat S(a))=0$, for all but a finite number of indices $j\in I_C$, see \cite[Definition 8.6]{Sim1}. Fix $a\in I_C$ and  let $0 \to \widehat S(a) \to \widehat E(a) \to \widehat E_1\to \dots $ be a minimal injective resolution of $\widehat S(a)$ in $C^{op}\textrm{-Comod}$, with $\widehat E_1 \cong  E(\soc(\widehat E (a)/\widehat  S(a)))$. Given $j\in I_C$, we denote by $\mu_1( \widehat  S(j), \widehat  S(a))$ the number of times the comodule $\widehat  E(j)$ appears as a direct summand in $\widehat E_1$. Since $C$ is assumed to be pointed,  $  \dim_K \End_C\widehat  S(j)=1$ and 
 $
  \mu_1(\widehat  S(j), \widehat  S(a))=  \dim_K \Ext^1_C(\widehat  S(j), \widehat  S(a)),   
  $ 
 by \cite[(4.23)]{Sim5}. Thus the injective $C^{op}$-comodule $\widehat E_1$ is socle finite, by the observation made earlier, and it follows that the simple right $C$-comodule $\widehat  S(a)$ is finitely copresented.  This shows that (b3) implies (b4). Since the inverse implication follows in a similar way, the proof of  (b) is complete.

  (c) Apply (b) and the   easily seen fact that   simple right comodules over any right locally artinian coalgebra  are finitely copresented. 
  \end{proof}    \medskip

Following   \cite{CKQ}, we get the following important result.\medskip
  
   {\sc Proposition 1.16.} \textit{ Let $C$ be a $K$-coalgebra and $D$ the duality functors} \textrm{(1.3)}.

  \textrm{(a)}   \textit{The transpose equivalence }   \textit{of Theorem}  1.8 (d), \textit{defines the equivalence} \vspace{-0.3cm}
  $$
   \textrm{Tr}_C  :  C\mbox{-}\overline{\textrm{Comod}}^\bullet_{fc}   \, \, \mapr{\simeq}{}\, \,  C^{op}\mbox{-}\overline{\textrm{comod}}_{fc} ,\vspace{-0.2cm}
   $$
     \textit{and  together with the duality} $D:  C^{op} \textrm{-comod}_{fc}  \, \, \mapr{\simeq}{}\, \,C \textrm{-comod}_{f\mathcal{P}}$  \textit{defined by} (1.3), \textit{induces the translate  operator}\vspace{-0.4cm}
     $$
 \tau^-_C= D \textrm{Tr}_C  :  C\mbox{-} {\textrm{Comod}}^\bullet_{fc}   \, \, \mapr{}{}\, \,  C\mbox{-}{\textrm{comod}}_{f\mathcal{P}} ,\vspace{-0.3cm}\leqno(1.17)
   $$
      \textit{and an equivalence of quotient categories}  $\overline{\tau}^-_C=D \textrm{Tr}_C:  C\mbox{-}\overline{\textrm{Comod}}^\bullet_{fc}   \, \, \mapr{\simeq}{}\, \,  C\mbox{-}\underline{\textrm{comod}}_{f\mathcal{P}}$.   
       \textit{Moreover,  for any $M$ in}  $\textrm{Comod}^\bullet_{fc}$, \textit{with a presentation}  (1.4), \textit{the following sequence}\vspace{-0.2cm}
   $$
  0  \, \mapr{}{}\,    (\widetilde D M)^+   \, \mapr{}{}\,  (\widetilde DE_1)^+\, \mapr{ (\widetilde Dg)^+}{}\,  (\widetilde D E_0)^+  \, \mapr{}{}\,  \tau^-_C M \, \mapr{}{}\,0\vspace{-0.2cm} \leqno(1.18)
  $$
 \textit{ is exact in} ${C^*}^{op}\mbox{-}\textrm{PC} $ \textit{and the comodule} $ \tau^-_C M  $ lies in $C \textrm{-comod}_{f\mathcal{P}}\equiv {C^*}^{op}\mbox{-}\textrm{rat}_{fp}\subseteq  {C^*}^{op}\mbox{-}\textrm{PC} $.

   \textrm{(b)}   \textit{The   duality} (1.3)  \textit{restricts to the duality} $D: C \textrm{-comod}_{f\mathcal{P}}  \, \, \mapr{\simeq}{}\, \,C^{op} \textrm{-comod}_{fc} $   \textit{and together with the  transpose operator }  $ \textrm{Tr}_{C^{op}}  :  C^{op}\mbox{-} {\textrm{comod}}_{fc}  \, \, \mapr{ }{}\, \, C\mbox{-} {\textrm{Comod}}^\bullet_{fc}  $   \textit{defines the translate operator} \vspace{-0.5cm}
  $$
  \tau_C= \textrm{Tr}_{C^{op}}D  : C \textrm{-comod}_{f\mathcal{P}}   \, \, \mapr{ }{}\, \,  C\mbox{-} {\textrm{Comod}}^\bullet_{fc}  ,\vspace{-0.4cm}\leqno(1.19)
   $$
      \textit{and induces the  equivalence of quotient categories} 
           $\overline{\tau}_C=  \textrm{Tr}_{C^{op}}D: C\mbox{-}\underline{\textrm{comod}}_{f\mathcal{P}}   \, \, \mapr{\simeq}{}\, \,  C\mbox{-}\overline{\textrm{Comod}}^\bullet_{fc} $ \textit{ that is quasi-inverse to the equivalence} $\overline{\tau}^-_C=D \textrm{Tr}_C:  C\mbox{-}\overline{\textrm{Comod}}^\bullet_{fc}   \, \, \mapr{\simeq}{}\, \,  C\mbox{-}\underline{\textrm{comod}}_{f\mathcal{P}}$ \textit{in} (a). 
  
\textrm{(c)}  \textit{Let   $ M$ be an indecomposable  comodule  in} $C \textrm{-Comod}^\bullet_{fc}$.  \textit{Then $\tau^-_CM=0$ if and only if $M$ is injective. If   $\tau^-_CM\neq 0$ then $\tau^-_CM $ is indecomposable, non-projective, of finite $K$-dimension, and there is an isomorphism $M \cong \tau_C\tau^-_CM$}. 

 \textrm{(d)}  \textit{Let   $ N$ be an indecomposable  comodule  in} $C \textrm{-comod}_{f\mathcal{P}}$. \textit{Then $\tau_CN  =0$ if and only if $N$ is projective. If   $\tau _CN\neq 0$ then $\tau _CN  $ is indecomposable, non-injective, finitely copresented,  and there is an isomorphism $N \cong \tau^-_C\tau_C N $.
 } \medskip
 
\begin{proof}  By Lemma 1.15(a),  the duality  $D:  C^{op} \textrm{-comod}   \, \, \mapr{\simeq}{}\, \,C \textrm{-comod} $   (1.3) restricts to the duality  $D:  C^{op} \textrm{-comod}_{fc}   \, \, \mapr{\simeq}{}\, \,C \textrm{-comod}_{f\mathcal{P}} $. One also shows, by applying foregoing definitions,  that a homomorphism $f:X\to X'$ in $C^{op} \textrm{-comod} _{fc}$ has a factorisation through a socle-finite injective comodule if and only if the homomorphism $D(f): D(X)\to D(X')$ in $C \textrm{-comod}_{f\mathcal{P}}$ belongs to $\mathcal{P}(D(X), D(X'))$. This shows that   the duality  $D:  C^{op} \textrm{-comod}_{fc}  \, \, \mapr{\simeq}{}\, \,C { \textrm{-comod}}_{f\mathcal{P}}$ induces an equivalence of quotient categories  $D:  C^{op} \overline{ \textrm{-comod}}_{fc}  \, \, \mapr{\simeq}{}\, \,C\underline{ \textrm{-comod}}_{f\mathcal{P}}$. It follows from the definition of the category $C\mbox{-}\overline{\textrm{Comod}}^\bullet_{fc}$ that the transpose equivalence     of Theorem   1.8 (d),  defines the equivalence 
   $
   \textrm{Tr}_C  :  C\mbox{-}\overline{\textrm{Comod}}^\bullet_{fc}   \, \, \mapr{\simeq}{}\, \,  C^{op}\mbox{-}\overline{\textrm{comod}}_{fc} .
   $ 
   This together with the earlier observation implies (a) and (b). 
   
   The statements (c) and (d) are obtained by a straightforward  calculation and by using the definition of translates $\tau_C$ and $\tau_C^-$, see  \cite{CKQ} and consult \cite{ARS}. The details are left to the reader. \smallskip
    \end{proof} 
        
    Following the terminology of representation theory of finite-dimensional algebras (see \cite{ASS}, \cite{ARS}, \cite{Si92}) we call the operators $ \tau_C= \textrm{Tr}_{C^{op}}D$ (1.19)  and $ \tau^-_C= D\textrm{Tr}_C$ (1.17), the \textbf{Auslander-Reiten translations} of $C$. It follows from  Theorem 1.16 that the image of $\tau^-_C$ is the subcategory $C \textrm{-comod}_{f\mathcal{P}}$ of the category  $C \textrm{-comod}$.
    
   By applying Theorem 4.2 and Corollary 4.3 in  \cite{CKQ}, we get the following  important  result on the existence of almost split sequences in the category $C \textrm{-Comod} _{fc}$ of (socle) finitely copresented left $C$-comodules, under some assumption on the Gabriel quiver   ${}_CQ$ of $C$. \medskip
    
     {\sc Theorem 1.20.}  \textit{Let $C$ be a $K$-coalgebra such that  its left Gabriel quiver}   ${}_CQ$  \textit{is left locally bounded, that is,  for every vertex $a$ of}   ${}_CQ$  \hbox{\textit{there is only a finite number of arrows $j\,\mapr{}{}\, a$ in} ${}_CQ$}.

\textrm{(a)}  \textit{The following  inclusion holds} $C \textrm{-comod} \subseteq C \textrm{-Comod} _{fc}$.

\textrm{(b)}   \textit{For any indecomposable non-injective comodule $M$ in} $C \textrm{-Comod}^\bullet_{fc}$,  \textit{there exists a unique almost split sequence}\vspace{-0.3cm}
$$  
0  \, \mapr{}{} \,\   M   \,\mapr{}{}  \,M' \,\mapr{}{} \,\tau^{-}_CM  \,
\mapr{}{} \, 0  \vspace{-0.2cm}  \leqno(1.21)
$$ 
\textit{in} $C\Comod_{fc}$, \textit{with a finite-dimensional  indecomposable comodule $\tau^{-}_CM$ lying in} $C \textrm{-comod}_{f\mathcal{P}}$. \textit{The sequence}  \textrm{(1.21)}  \textit{is almost split in the whole comodule category} $C\Comod$. 

\textrm{(c)}   \textit{For any   indecomposable non-projective comodule $N$ in} $C \textrm{-comod}_{f\mathcal{P}}\subseteq C \textrm{-Comod} _{fc}$,  \textit{there exists a unique almost split sequence}\vspace{-0.2cm}
$$  
0  \, \mapr{}{} \,\  \tau_C N   \,\mapr{}{}  \,N' \,\mapr{}{} \,N  \,
\mapr{}{} \, 0  \vspace{-0.2cm}   \leqno(1.22)
$$ 
\textit{in} $C\Comod_{fc}$, \textit{with an indecomposable comodule $\tau _CN$ lying in} $C \textrm{-Comod}^\bullet_{fc}$.  \textit{The sequence}  \textrm{(1.22)}  \textit{is almost split in the whole comodule category} $C\Comod$. 

\textrm{(d)}  \textit{If,  in addition, $C$ is left semiperfect  then} $C\Comod^\bullet_{fc}=  C\Comod_{fc}$, $C \textrm{-comod}_{f\mathcal{P}}= C \textrm{-comod}$,  \textit{the Auslander-Reiten translate operators act as follows} \vspace{-0.25cm}
$$
  C\mbox{-} {\textrm{Comod}}_{fc}   \, \,\RightLeftarr{\tau^-_C}{\tau_C}\, \,  C\mbox{-}{\textrm{comod}} \vspace{-0.25cm}
$$
\textit{and the almost split sequences} (1.21) \textit{and} (1.22) \textit{do  exist in the category} $C\Comod_{fc}$,  \textit{for any indecomposable non-injective comodule $M$ in }  $C\Comod_{fc}$ \textit{and for any indecomposable non-projective comodule $N$ in} $C \textrm{-comod}$. \textit{Moreover,  if the comodule $M$ lies in} $C \textrm{-comod}$ \textit{then the  almost split sequence} (1.21 ) \textit{lies in } $C \textrm{-comod}$. 

    \medskip
 
\begin{proof}  (a) As in the proof of Lemma 1.15 (b), we conclude from  the assumption that ${}_CQ$   is left locally bounded that 
every simple left $C$-comodule admits a minimal socle-finite injective copresentation (1.4). Hence (a) follows as in Lemma 1.15 (b).

The statements (b) and (c) follow from Theorem 4.2 and Corollary 4.3 in  \cite{CKQ}, because any comodule $M$  lying in   $C\Comod_{fc}$ is quasi-finite, Proposition 1.16 (a) yields that $\tau^-_CM$ lies in $C \textrm{-comod}_{f\mathcal{P}}$, for any indecomposable comodule $M$ in $C\textrm{-Comod}^\bullet _{fc}$,  and the following inclusions hold   $C \textrm{-comod}_{f\mathcal{P}} \subseteq C \textrm{-comod} \subseteq C \textrm{-Comod} _{fc}$, by (a).

(d)  Assume that $C$ is left semiperfect and  let $M$ be an indecomposable comodule in $C\Comod_{fc}$ with a minimal socle-finite injective copresentation (1.4).  By Theorem 1.8,   the induced sequence (1.10) is exact and the comodules $ \nabla_C(E_0)$ and $ \nabla_C(E_1)$ lie in $C^{op}\inj$. Since $C$ is left semiperfect, the comodules $ \nabla_C(E_0)$ and $ \nabla_C(E_1)$ are finite-dimensional and, hence, $\dim_K   \textrm{Tr}_C  (M)  $ is finite, for any comodule $M$ in $C\Comod_{fc}$. It follows that $C\Comod^\bullet_{fc}=  C\Comod_{fc}$.  Since $C$ is left semiperfect, any comodule $N$ in $ C\comod$ has a projective presentation  $P_1\, \mapr{}{}\, P_0 \, \mapr{}{}\, N \, \mapr{}{}\, 0$, with $P_1, P_0$ finite-dimensional projective $C$-comodules.  It follows that $N$ lies in $C \textrm{-comod}_{f\mathcal{P}}$ and, hence,  the equality
$C \textrm{-comod}_{f\mathcal{P}}= C \textrm{-comod}$ holds.  This finishes the proof of the theorem.
\end{proof}
 \smallskip
 
    {\sc Corollary 1.23.} \textit{Let $C$ be a pointed   $K$-coalgebra such that the  left Gabriel quiver}   ${}_CQ$  \textit{of $C$ is both  left and right  locally bounded}.  
    
 \textrm{(a)}   \textit{ The inclusions} $C\mbox{-} {\textrm{comod}}_{f\mathcal{P}} = C\mbox{-} {\textrm{comod}}\subseteq C\mbox{-} {\textrm{Comod}}_{fc}$ \textit{hold  and the Auslander-Reiten translate operators act as follows} \vspace{-0.3cm}
$$
  C\mbox{-} {\textrm{Comod}}^\bullet_{fc}   \, \,\RightLeftarr{\tau^-_C}{\tau_C}\, \,  C\mbox{-}{\textrm{comod}}.  \vspace{-0.3cm}
$$

\textrm{(b)}  \textit{For any indecomposable non-injective comodule $M$ in} $C \textrm{-Comod}^\bullet_{fc}$,  \textit{there exists a unique almost split sequence}\vspace{-0.3cm}
$$  
0  \, \mapr{}{} \,\   M   \,\mapr{}{}  \,M' \,\mapr{}{} \,\tau^{-}_CM  \,
\mapr{}{} \, 0  \vspace{-0.2cm}  
$$ 
\textit{in} $C\Comod_{fc}$, \textit{with an   indecomposable comodule $\tau^{-}_CM$ lying in} $C \textrm{-comod}$. 

 \textrm{(c)}       \textit{For any   indecomposable non-projective comodule $N$ in} $C \textrm{-comod}$,  \textit{there exists a unique almost split sequence}\vspace{-0.2cm}
$$  
0  \, \mapr{}{} \,\  \tau_C N   \,\mapr{}{}  \,N' \,\mapr{}{} \,N  \,
\mapr{}{} \, 0  \vspace{-0.2cm}  
$$ 
\textit{in} $C\Comod_{fc}$, \textit{with an indecomposable comodule $\tau _CN$ lying in} $C \textrm{-Comod}^\bullet_{fc}$.

(d)  \textit{The exact sequences in}  \textrm{(b)}    \textit{ and }  \textrm{(c)}  \textit{are  almost split in the whole comodule category} $C\Comod$. 
  \medskip
  
 \begin{proof}    
  Apply Lemma 1.15 and Theorem 1.20.
  \end{proof}    \medskip
  
      {\sc Remark 1.24.}   Under the assumption that the left  Gabriel quiver ${}_CQ$  of $C$ is both  left and right  locally bounded the almost split sequences  (1.21) and (1.22) lie in $C\Comod_{fc}$. If we drop the assumption then the term $\tau_CM $ lies in  $C\mbox{-} {\textrm{comod}}_{f\mathcal{P}} \subseteq  C\mbox{-} {\textrm{comod}}$, but not necessarily lies in $C\Comod_{fc}$.  Since the category $C\Comod_{fc}$ is the most important part of $C\Comod$ in the study of the tameness of $C$ (see \cite{Sim08}),  we are mainly interested in the existence of  almost split sequences in $C\Comod_{fc}$.       
       \medskip

  \vfill\eject 
 
  Now we illustrate the existence of almost split sequences discussed in  Corollary 1.23 by the following example.  \medskip
  
  \def\vedge#1{{\buildrel{#1} \over {\hbox to
20pt{\hspace{-0.2em}$-$\hspace{-0.2em}$-$\hspace{-0.2em}$-$ }}}}
 
    {\sc Example 1.25.}  Let   $Q= (Q_0, Q_1)$ be the  infinite locally Dynkin quiver   \vspace{-0.3cm}
$$
     \,\,  Q:   
\begin{array}{cccccccccc}
&& \hspace{-0.6em}{\scriptstyle {-1}}  &   &    &    &    &   & &\vspace{-1.5ex}\\
&& \bullet  &   &    &    &    &   & &\vspace{-0.8ex}\\
&& \hspace{-0.0em}\mapupp{}{} & &  &  &  & & &\vspace{-0.6ex}\\
\bullet  &  \hspace{-1em}\mapl{}{}  \hspace{-1em}&\circ   &
\hspace{-1em}\mapr{}{}  \hspace{-1em} &\bullet   &\hspace{-0.81em}\mapr{}{}  \hspace{-0.0em} \,\,
 \cdots\, \, \hspace{-0.01em}\mapr{}{}  \hspace{-1.1em}& \bullet   & \hspace{-1.3em}\mapr{}{}  \hspace{-1em} & \bullet  &\hspace{-0.7em}\mapr{}{}    \hspace{-0.5em}\,\,\bullet \,\, \hspace{-0.5em}\mapr{}{}   \hspace{-0.1em}\,\,
\ldots \, 
\vspace{-1.4ex}\\
{\scriptstyle { 0}}  &  &{\scriptstyle { 1}}  &
\hspace{-1cm}
  & {\scriptstyle {2}}  & \hspace{-1cm} \,\,
  \, \, \hspace{-2cm}& {\scriptstyle {{s-1}}}  &\hspace{-1cm}& {\scriptstyle { s}}  &\hspace{-1cm}{\scriptstyle {{ s+1}}}\, \end{array}  \vspace{-0.3cm}
$$ 
  of type $\mathbb{D}_\infty$ and let $C=   K^\SQR Q$ be the hereditary path $K$-coalgebra of $Q$.   Then  $Q_0 =   \{-1, 0,1,2,\ldots\}$ and $C$ has the   $Q_0 \times Q_0 $ matrix form 
 $$
C =\left[\scriptsize{\begin{matrix}
  &K&0 &0 & 0 & 0& 0& 0&0&0&\ldots\cr
   &0&K &0 & 0 & 0& 0& 0&0&0&\ldots\cr
  &K&K&K & K & K& K& K&K&K&\ldots\cr
  &0&0&0 & K & K& K& K&K&K&\ldots\cr
  &0&0&0 & 0 &  K& K& K&K&K&\ldots\cr
  &0&0&0 & 0& 0& K& K&K&K&\ldots\cr
  &0&0&0& 0& 0& 0& K&K&K&\ldots\cr
  &0&0&0&0& 0& 0& 0& K&K&\ldots\cr
  &0&0&0&0& 0& 0& 0& 0&K&\ldots\cr
 &\vdots & \vdots & \vdots & \vdots& \vdots& \vdots&\vdots &\vdots
&\vdots   &\ddots &\cr
\end{matrix}}
\right]
$$ 
and consists of the   triangular $Q_0 \times Q_0 $ square matrices with coefficients in $K$  and   at most  finitely many 
non-zero entries.  Then $\soc{}_CC =   \bigoplus _{j\in Q_0} S(j)$, where $S(n) = Ke_n$ is the simple subcoalgebra of $C$ spanned by the matrix $e_n\in C$ with $1$ in the $ n\times n$ entry, and zeros elsewere. Note that $e_n$ is a group-like element of $C$.

  Since the left Gabriel quiver ${}_CQ$ of $C$ is the quiver $Q$, it follows from  Lemma 1.15 (b) that   every simple right $C$-comodule is finitely copresented and the statements (a) and (c) of Corollary 1.23 hold for $C=   K^\SQR Q$.
  
  The coalgebra $C$ is right locally artinian, right semiperfect, representation-directed in the sense of \cite{Sim4},  and left pure semisimple, that is, every left $C$-comodule is a direct sum of finite-dimensional comodules (see \cite{NowSim}, \cite{Sim1},  and \cite{Sim2}).  It follows that $C\Comod^\bullet_{fc} = C\comod^\bullet$\linebreak and every indecomposable non-projective comodule $N$ in $C\comod$ admits an almost split sequence  $0  \, \mapr{}{} \,\  \tau_C N   \,\mapr{}{}  \,N' \,\mapr{}{} \,N  \,
\mapr{}{} \, 0 $ 
in $C\comod$. 
  
  Under the identification $C\comod= \rep_K(Q)$ of left $C$-comodules and $K$-linear representations of the quiver $Q$ (see \cite{CKQ}, \cite{Sim1}, \cite[(3.1)]{Sim7}), the Auslander-Reiten translation quiver $\Gamma(C\comod)$  of  $C\comod$ has   four connected components (two of them are finite and two are infinite),  and  $\Gamma(C\comod)$ has the following form  (see \cite{NowSim} and \cite{Sim2})
    \vfill\eject 

\begin{center}
\underline{{\sc Figure 0.} {\sc The Auslander-Reiten quiver of 
the category $ C\comod \cong  \rep_K( Q))$ } }
\end{center}

  $$
\hspace{-0.3cm}\scriptsize{
\begin{array}{ccccccccccccccccccccccccccccccccccc}
  \shp \hspace{-15.5cm}\hidewidth_{0}{\Bbb I}_{0}\hidewidth  &   &&  & &  &&  & & && &  &    && & & &   && & &    && &
&  && &  &    &&
 \vspace{2ex} \\
 &  &   &&\shk\hspace{-17.5cm} \hidewidth _{-1}{\Bbb I}_{-1} \hidewidth  & &  &&  & & && &  &    && & &
&   && & &    && & &  && &  &    && 
 \vspace{0.4cm}\\
  &  &   &&  
  & \ldots &
   \hidewidth\hbox{ - - - }\hidewidth&\hidewidth _{0}{\Bbb I}_6 \hidewidth& & \hidewidth\hbox{ - - - - - - - - }\hidewidth &     
 &\hidewidth _{-1}{\Bbb I}_5 \hidewidth& &\hidewidth\hbox{ - - - - - - - - }\hidewidth& &
\hidewidth_{0}{\Bbb I}_4 \hidewidth &&
\hidewidth\hbox{ - - - - - - - - }\hidewidth & &\hidewidth _{-1}{\Bbb I}_3 \hidewidth &&
\hidewidth\hbox{ - - - - - - - - }\hidewidth & & \hidewidth_{0}{\Bbb I}_2 \hidewidth &&
\hidewidth\hbox{ - - - - - - - - }\hidewidth&  & \hidewidth _{-1}{\Bbb I}_1 \hidewidth  && & &  && 
\medskip \\
    &  & & &    &  &   \nearrow
&  &\searrow&  &   \nearrow &  &
\searrow & & \nearrow& &\searrow&   &\nearrow &    &\searrow &  &\nearrow &   
&\searrow &&\nearrow &    &\searrow & & &  && 
\medskip \\ 
 &    & &&
&& &\ldots&  \to &\hidewidth_{4}{\Bbb K}_5\hidewidth&\to & \hidewidth_{0}{\Bbb I}_5\hidewidth& 
\to&\hidewidth _{4}{\Bbb K}_5 \hidewidth&\to&\hidewidth_{-1}{\Bbb I}_4\hidewidth&\to&
\hidewidth_{3}{\Bbb K}_4 \hidewidth&\to &\hidewidth_{0}{\Bbb I}_3\hidewidth&\to&\hidewidth _{2}{\Bbb K}_3\hidewidth
&\to& 
\hidewidth_{-1}{\Bbb I}_2 \hidewidth
&\to &\hidewidth _{1}{\Bbb K}_2 \hidewidth &\to &\hidewidth _{0}{\Bbb I}_1\hidewidth   & \to
&\hidewidth_1{\Bbb I}_1\hidewidth & &  && 
\medskip \\
        &  &  &  &    &  && 
&  \nearrow  &  &\searrow&  &\nearrow
 &  &\searrow&   &\nearrow &    &\searrow &  &\nearrow &    &\searrow &&\nearrow
&    &\searrow && \nearrow&& &  && 
\medskip \\ 
     & &  &&  & &  && & \hidewidth\hbox{ - - - - - - - - }\hidewidth 
&    &\hidewidth_{4}{\Bbb K}_6 \hidewidth&
   &\hidewidth\hbox{ - - - - - - - - }\hidewidth& &\hidewidth _{3}{\Bbb K}_5\hidewidth &&
\hidewidth\hbox{ - - - - - - - - }\hidewidth  & &\hidewidth _{2}{\Bbb K}_4 \hidewidth &&
\hidewidth\hbox{ - - - - - - - - }\hidewidth & &\hidewidth_{1}{\Bbb K}_3 \hidewidth &&
\hidewidth\hbox{ - - - - - - - - }\hidewidth&  & \hidewidth _{1}{\Bbb I}_2 \hidewidth &&& &  && 
\medskip \\
    &  &    &  && 
 &    & &&  &\nearrow& 
&   \searrow
 & &\nearrow&    & \searrow &  &\nearrow &    &\searrow &&\nearrow
 &    &\searrow&& \nearrow&& 
&   & &  && 
  \medskip\\
    & &  &&  
&&   &&   &
&&\hidewidth\hbox{ - - - - - - - - }\hidewidth &
   &\hidewidth   _{3}{\Bbb K}_6 \hidewidth& & \hidewidth\hbox{ - - - - - - - - }\hidewidth &&\hidewidth
_{2}{\Bbb K}_5\hidewidth  && \hidewidth\hbox{ - - - - - - - - }\hidewidth &&\hidewidth _{1}{\Bbb
K}_4 \hidewidth &&  \hidewidth\hbox{ - - - - - - - - }\hidewidth & &\hidewidth_{1}{\Bbb I}_3 \hidewidth& &&&
& &  && 
\medskip \\
           & &&  && 
&    & &&   &   &  &  
    \nearrow   & &\searrow &&\nearrow
 &    &\searrow&& \nearrow&&\searrow
&    & \nearrow&&   &   &  &  & &  && 
 \medskip  \\
&&   && &&
  &&   & && &
   & \hidewidth\hbox{ - - - - - - - - }\hidewidth &&\hidewidth  _{2}{\Bbb K}_6 \hidewidth& &
\hidewidth\hbox{ - - - - - - - - }\hidewidth &&\hidewidth _{1}{\Bbb K}_5 \hidewidth &&
\hidewidth\hbox{ - - - - - - - - }\hidewidth&&\hidewidth _{1}{\Bbb I}_4 \hidewidth &&   & && &
& &  && 
  \medskip \\
          &  & & 
&   &  &&    &   &  &&&
        & &\nearrow
 &    &\searrow& & \nearrow& &\searrow
&   & \nearrow &&    &   &  &&&& &  && 
\\
 &  & & 
&     & &&    &  &    &  &  & 
 && 
 &   \hidewidth\hbox{ - - - - - - - - }\hidewidth  & &\hidewidth _{1}{\Bbb K}_6\hidewidth&  &
\hidewidth\hbox{ - - - - - - - - }\hidewidth &  & \hidewidth   _{1}{\Bbb I}_5\hidewidth & &&    & 
&    &  &  & & &  && 

\\
          & &&  && 
&    & &&   &   &  &  
        & &  &&\nearrow
 &    &\searrow&& \nearrow&& 
&    &  &&   &   &  &  & &  && 
  \medskip  \\
&&   && &&
  &&   & && &
     &   &&   & &\hidewidth\hbox{ - - - - - - - - }\hidewidth 
  && \hidewidth_{1}{\Bbb I}_6 \hidewidth &&
 &&  &&   & && && &  && 
\\
           & &&  && 
&    & &&   &   &  &  
        & &  && 
 &    &\hidewidth\nearrow\hidewidth&& && 
&    &  &&   &   &  &  & &  && 
\vspace{-0.5cm}\\
     & 
 &    & &\ddots&  & &  
&   \ddots  &  &&     &  &    &  &  &
 && 
 &    & & &  & & 
&     & &&    &  &    &  &  &
\medskip \\
 & \hidewidth  \searrow\hidewidth&&\hidewidth\nearrow\hidewidth&    & \searrow &  &\nearrow &    &\searrow &&
 &    & 
 &  &  &    &  &&
 &    && & &  &&&&&&&&
&  \medskip\\
 \ldots      &\hidewidth\hbox{ - - - }\hidewidth& \hidewidth _{3}{\Bbb I}_{7}\hidewidth& & 
     \hidewidth\hbox{ - - - - - - - - }\hidewidth && \hidewidth _{2}{\Bbb I}_6\hidewidth & &
\hidewidth\hbox{ - - - - - - - - }\hidewidth && 
\hidewidth{\Bbb I}_5\hidewidth  &&   & &&   &&&&
 && &&
   &&   & &  & & &&
\medskip\\
       &\nearrow
 &    &\searrow& & \nearrow
 &    &\searrow& & \nearrow& &\searrow
&   &    &   &  &&&&& & 
&   &    &   &  &&&   & &
  &   & &
\medskip\\
 \ldots    & & \hidewidth\hbox{ - - - - - - - - }\hidewidth &&\hidewidth _{3}{\Bbb
I}_{6} \hidewidth && \hidewidth\hbox{ - - - - - - - - }\hidewidth &&\hidewidth _{2}{\Bbb I}_{5} \hidewidth
 & & \hidewidth\hbox{ - - - - - - - - }\hidewidth &&
\hidewidth {\Bbb I}_4  \hidewidth&&   && &
&&&&&&&&&&&& &&&&
\medskip \\
 &   \searrow&&\nearrow&    & \searrow &  &\nearrow &    &\searrow &&
\nearrow&    & 
\searrow &  &  &    &  &&
 &    && & &  &&&&&&&&
&  \medskip\\
 \ldots      &\hidewidth\hbox{ - - - }\hidewidth&\hidewidth  _{4}{\Bbb I}_{6} \hidewidth& &
\hidewidth\hbox{ - - - - - - - - }\hidewidth &&\hidewidth _{3}{\Bbb I}_{5}\hidewidth& & 
     \hidewidth\hbox{ - - - - - - - - }\hidewidth &&  \hidewidth_{2}{\Bbb I}_4 \hidewidth& &
\hidewidth\hbox{ - - - - - - - - }\hidewidth && 
\hidewidth{\Bbb I}_3 \hidewidth &&   & &&   &&&&
 && &&
   &&   & &
\medskip \\
         & \nearrow& &\searrow &&\nearrow
 &    &\searrow&&  \nearrow& &\searrow &&\nearrow
 &    &\searrow&& && 
&    &   &   &  & & & 
&    &   &   &  & 
 & &
   \medskip  \\
 \ldots    & & \hidewidth\hbox{ - - - - - - - - }\hidewidth &&\hidewidth _{4}{\Bbb
I}_{5} \hidewidth && \hidewidth\hbox{ - - - - - - - - }\hidewidth &&\hidewidth _{3}{\Bbb I}_{4}\hidewidth 
 & & \hidewidth\hbox{ - - - - - - - - }\hidewidth &&
\hidewidth_{2}{\Bbb I}_3 \hidewidth && \hidewidth\hbox{ - - - - - - - - }\hidewidth &&\hidewidth{\Bbb I}_2 \hidewidth&
&&&&&&&&&&&& &&&&
\medskip \\
 &   \searrow&&\nearrow&    & \searrow &  &\nearrow &    &\searrow &&
\nearrow&    & 
\searrow &  &\nearrow &    &\searrow &&
 &    && & &  &&&&&&&&
&  
   \medskip \\
 \ldots&   \hidewidth\hbox{ - - - }\hidewidth&\hidewidth
_{5}{\Bbb I}_{5} \hidewidth
    &&\hidewidth\hbox{ - - - - - - - - }\hidewidth &  & \hidewidth _{4}{\Bbb I}_4\hidewidth  & 
&\hidewidth\hbox{ - - - - - - - - }\hidewidth& & \hidewidth_{3}{\Bbb I}_3  \hidewidth&&
\hidewidth\hbox{ - - - - - - - - }\hidewidth  & &\hidewidth _{2}{\Bbb I}_2 \hidewidth &&
\hidewidth\hbox{ - - - - - - - - }\hidewidth&  & \hidewidth {\Bbb I}_1 \hidewidth  &&& &&&&&&&&&&&&
\end{array}
}
$$
Here we use the terminology and notation introduced in \cite[pp. 470--472]{NowSim} and   \cite[Section 6]{Sim2}.   Recall that the vertices of the Auslander-Reiten translation quiver $\Gamma(C\comod)$   are representatives of  the indecomposable left $C$-comodules in $C\comod$ and the existence of an arrow $X\to Y$ in   $\Gamma(C\comod)$ means that there exists an  irreducible morphism  $f:X\to Y$ in $C\comod$, see also  \cite{Sim6}.

  Each of the two finite components of $\Gamma(C\comod)$  contains  precisely one indecomposable simple projective $C$-comodule; namely the comodule $_{0}{\Bbb I}_{0}$ and $_{-1}{\Bbb I}_{-1}$, respectively.   Each of the two infinite components  
 contains    no non-zero projective objects.   

The indecomposable injective left $C$-comodules  form the right hand   section  \vspace{-0.2cm}
 $$
 \begin{array}{r}
 {}_0{\Bbb I}_1\,\,\quad\qquad  \\
  \downarrow{}{}\,\,\quad\qquad  \\
 \ldots \to {} _1{{\Bbb I}}_6\to  {}_1{{\Bbb I}}_5\to  {}_1{{\Bbb I}}_4\to  {}_1{{\Bbb I}}_3 \to  {}_1{{\Bbb I}}_2\to  {}_1{{\Bbb I}}_1  \leftarrow {}_{-1}{\Bbb I}_1 
 \end{array}\vspace{-0.2cm}\leqno(*)
 $$  
 of the infinite upper component of $\Gamma(C\comod)$, and  the indecomposable left $C$-comodules $M$ in $ C\comod$ such that $\dim_K \mathrm{Tr}_C M$ is infinite  are the two simple projective comodules $ {}_0{\Bbb I}_0$, $ {}_{-1}{\Bbb I}_{-1}$ and the comodules   lying on the  infinite   right hand   section  \vspace{-0.2cm}
 $$ 
 \ldots \to  {{\Bbb I}}_6\to {{\Bbb I}}_5\to  {{\Bbb I}}_4\to   {{\Bbb I}}_3 \to  {{\Bbb I}}_2\to   {{\Bbb I}}_1 
 \vspace{-0.2cm}\leqno(**)
 $$  
 of the infinite  lower component of $\Gamma(C\comod)$.  It follows that an indecomposable left $C$-comodule $M$ lies in the category $C\Comod^\bullet_{fc} = C\comod^\bullet$ if and only if $M$ lies in the infinite upper component of $\Gamma(C\comod)$ or $M$ lies in the infinite lower component of $\Gamma(C\comod)$, but does not lie on the infinite section $(**)$.  Every indecomposable  comodule $N$ lying in one of the infinite components of $\Gamma(C\comod)$ has an almost split sequence $0  \, \mapr{}{} \,\  \tau_C N   \,\mapr{}{}  \,N' \,\mapr{}{} \,N  \,
\mapr{}{} \, 0 $ 
in $C\comod$ and it is given by  the mesh $\Gamma(C\comod)$ terminating at $N$, compare with the examples given in Section 4.

   
 \section{Cartan matrix of an Euler  coalgebra and its inverses}

$\quad$ Throughout we assume that $K$ is an arbitrary   field and $C$ is a \textbf{pointed   $K$-coalgebra}. It follows that $C$ is basic and   there exists a direct sum decomposition $\soc\,{}_CC = \bigoplus _{j\in I_C} S(j)$ of the left socle $\soc\,{}_CC $ of $C$, where $I_C$ is an index set and $\{S(j)\}_{j\in I_C}$ is a set of pairwise non-isomorphic simple left $C$-coideals, see \cite{Ch}, \cite{CMo}, \cite{Sim0}. Then  $\{S(j)\}_{j\in I_C}$ is a set of representatives of  the isomorphism
classes of simple left $C$-comodules and $\dim_K S(j) = \dim_K \End_CS(j) = 1$, for  any $j\in I_C$, see \cite{Sim5}. 

For every
$j\in I_C$, let $E(j)=E(S(j))$  denote the  injective envelope of $S(j)$. It follows that $E(j)$ is indecomposable, ${}_CC = \bigoplus _{j\in I_C} E(j)$, and  there is a primitive idempotent
$e_{j}\in C^*$ such that  $E(j)\cong  e_jC$.
Working with right $C$-comodules, we have the  simple right $C$-comodules $ \widehat{S}(j)=DS(j)$,  with injective envelopes $\widehat{E}(j) =  \nabla_C(E(j))\cong Ce_j $ in $C^{op}\inj$, see Theorem 1.8.  Throughout we fix a set  $\{e_j\}_{j\in I_C}$ of primitive idempotents of $C^*$ such that $E(j)\cong  e_jC$, for all $j\in I_C$.

Following the representation theory of  finite-dimensional algebras, given a left $C$-comodule $M$ (viewed as a rational right $C^*$-module), we define its \textbf{dimension vector}\vspace{-0.2cm}
$$
\bdim\, M = [\dim_K Me_j]_{j\in I_C},\vspace{-0.2cm} \leqno(2.1)
$$
where $\dim_K Me_j$  has values in
$\mathbb{N}\cup\{\infty\}.$ Since $C$ is pointed,   $\dim_K Me_j =  \dim_K \Hom_C(M, e_jC) = \dim_K \Hom_C(M, E(j)) $ and the dimension vector  $\bdim\, M $ coincides with the composition length vector $\len M = [\ell_j(M)]_{j\in I_C}$ of $M$ (introduced in  \cite{Sim4}), where \vspace{-0.2cm}
$$
\ell_j(M) = \dim_K \Hom_C(M, E(j))= \dim_K Me_j  .\vspace{-0.2cm} \leqno(2.2)
$$
   It follows from \cite[Proposition 2.6]{Sim4} that  $\ell_j(M)= \dim_K Me_j  $ is the multiplicity the simple comodule $S(j)$ appears as a composition factor in the socle filtration  $
\soc^0M\subseteq \soc^1M\subseteq
\ldots 
\subseteq \soc^mM\subseteq
\ldots    
$  of $M$.   Following  \cite{Sim4}, $M$ is said to be \textbf{computable} if  the composition length  multiplicity $\ell_j(M)= \dim_K Me_j  $ of $S(j)$ in $M$ is finite, for every $j\in I_C$, or equivalently, $\bdim\, M\in \mathbb{Z} ^{I_C}$ (the product of $I_C$ copies of the infinite cyclic group $\mathbb{Z}$). A  pointed   coalgebra $C$ is defined to be   \textbf{computable} if the injective comodule $E(i)$ is computable, or equivalently, if  the dimension vector \vspace{-0.2cm}
$$
\textbf{e}(i)= \bdim\, E(i) = [\dim_K e_iCe_j ] _{j\in I_C}= [\dim_K \Hom_C(E(i), E(j)) ] _{j\in I_C} \vspace{-0.2cm}\leqno(2.3)
$$
has finite coordinates, 
 for every $i\in I_C$.  Note that the class of computable coalgebras contains left semiperfect coalgebras, right semiperfect coalgebras and the incidence coalgebras $ K^\SQR\! I$  of intervally finite posets $I$, see \cite{Sim4}. Moreover, if $C$ is computable and left cocoherent then the $K$-category
 $C\Comod_{fc}$ is  abelian, has enough injective objects, and is   $\Ext$-\textbf{finite}, that is, $\dim_K \Ext^m_C(M,N)$ is finite, for all $m\geq 0$ and all comodules $M$, $N$ in $C\Comod_{fc}$
  
Given a pointed   computable coalgebra $C$, with a fixed decomposition $\soc\,{}_CC = \bigoplus _{j\in I_C} S(j)$,  we define the left  \textbf{Cartan matrix} of $C$  to be the integral $  {I_C}\times  {I_C}$
matrix \vspace{-0.4cm}
$$
 \Car_C =[ \mathbf{c}_{ij}]_{i,j\in I_C} = 
  \left[\scriptsize{\begin{matrix}\cr\vspace{-0.7cm}\cr \vdots\cr \textbf{e}(i)\vspace{-0.2cm}\cr \vdots\end{matrix}}
 \right] 
 \in \mathbb{M}_{I_C}(\mathbb{Z}), \vspace{-0.2cm}\leqno(2.4)
$$
whose $i\times j$ entry   is the composition length multiplicity $\mathbf{c}_{ij}= \textbf{e}(i)_j = \dim_K e_iCe_j $
of
$S(j)$ in $E(i)$. In other words, the $i$th row of $ {\Large{\mbox{$\mathfrak{c}$}}}_C$ is the dimension vector $\textbf{e}(i)= \bdim\, E(i) $ of $E(i)$, see \cite[Definition 4.1]{Sim4}.
We say that a row (or a column) of a matrix is finite, if  the  number of its non-zero coordinates is  finite. A matrix is called \textbf{row-finite} (or \textbf{column-finite}) if each of its row (column) is finite.

\medskip

We start with the following simple observations.\medskip

  {\sc Lemma 2.5.} \textit{Let $C$ be a pointed   computable $K$-coalgebra, with a fixed decomposition $\soc\,{}_CC = \bigoplus _{j\in I_C} S(j)$, and let} $K_0(C) = K_0(C\comod)$ \textit{be  the Grothendieck group of} $C\comod$.
   
  \textrm{(a)}   \textit{Given a $C$-comodule $M$ in $C\Comod$, the dimension vector $\dim\, M$ has only a finite   number of   non-zero coordinates if and only if $\dim_K M$ is finite. If  $\dim_K M<\infty$, then $\dim_K M = \sum_{j\in I_C} \dim_K Me_j$.  }

   \textrm{(b)}  The map   $M\mapsto \bdim\, M$ \textit{ is an additive function on  short exact sequences in } $C\Comod$ \textit{and induces the group isomorphism} $\bdim: K_0(C)\,\,\mapr{\simeq}{} \mathbb{Z}^{(I_C)}$, $[M]\mapsto \bdim\, M$, \textit{where   $\mathbb{Z}^{(I_C)}$ is the direct sum of $I_C$ copies of $\mathbb{Z} $. The group $K_0(C)$ is free abelian with the basis $\{[S(j)]\}_{j\in I_C}$ corresponding via} $\bdim$ \textit{to the standard basis vectors} $e_j = \bdim\, S(j)$ \textit{of $\mathbb{Z}^{(I_C)}$.} 
  \medskip
 
\begin{proof} (a)  To prove the sufficiency, assume that  $\dim_K M$ is finite.  Then 

$M^* \cong \Hom_C(M, C) \cong \Hom_C(M,  \bigoplus _{j\in I_C} E(j)) \cong    \bigoplus _{j\in I_C} \Hom_C(M,  E(j)) \cong     \bigoplus _{j\in I_C}  Me_j$.

\noindent It follows that $\dim_K M=  \dim_K M^* =  \sum _{j\in I_C}  \dim_KMe_j$. Hence,   the sum is finite and $\bdim\, M$ has only a finite   number of   non-zero coordinates. The converse implication follows in  a similar way. 

(b) Since (2.2) yields $\bdim= \len$,   \cite{Sim1} and \cite{Sim4} apply.
\end{proof}
\medskip

   {\sc Lemma 2.6.} \textit{Let $C$ be a pointed   computable $K$-coalgebra and let $\Car_C \in \mathbb{M}_{I_C}(\mathbb{Z})$ be  the left Cartan matrix} \textrm{(2.4)}  \textit{of $C$}.
   
  \textrm{(a)}   \textit{$\Car_{C^{op}} =  \Car^{tr}_C$, that is, $\Car_{C^{op}} \in \mathbb{M}_{I_C}(\mathbb{Z})$ is  the transpose  of the matrix $\Car_C$. }

   \textrm{(b)}   \textit{The   $i$th row of the matrix  $\Car_C $ is finite if and only if  the indecomposable injective left $C$-comodule $E(i)$ is finite-dimensional.} 
   
    \textrm{(c)}   \textit{The   $j$th column  of the matrix  $\Car_C $ is finite if and only if  the indecomposable injective  right $C$-comodule $\widehat{E}(j)=  \nabla_C(E(j))$ is finite-dimensional.} 
  
\textrm{(d)}  \textit{The left Cartan matrix $\Car_C$ of $C$ is row-finite if and only if $C$ is right semiperfect}. 

\textrm{(e)}  \textit{The left Cartan matrix $\Car_C$ of $C$ is column-finite if and only if $C$ is left semiperfect}.\medskip
 
\begin{proof} (a) Let  $\Car_{C^{op}} =[ \overline{\mathbf{c}}_{ij}]_{i,j\in I_C} , \Car_{C} =[ {\mathbf{c}}_{ij}]_{i,j\in I_C}   \in \mathbb{M}_{I_C}(\mathbb{Z})$ be the Cartan matrices  (2.4), of  the coalgebra $C^{op}$ and $C$, respectively. We recall from Theorem 1.8 that there is a duality $ \nabla_C: C\inj \to C^{op}\inj$ and $\widehat{E}(j) =  \nabla_C(E(j))$. Hence, by applying (2.3),  we get 
$ \overline{\mathbf{c}}_{ij} = \dim_K\Hom_{C^{op}}( \widehat{E}(i), \widehat{E}(j)) =\dim_K\Hom_ C ( {E}(j),  {E}(i))  = \mathbf{c}_{ji}$, for every pair of elements $i,j\in I_C$. This yields the equality $\Car_{C^{op}} =  \Car^{tr}_C$.

(b)  We recall from (2.4) that the $i$th row of $\Car_C$ is the dimension vector $\textbf{e}(i) = \bdim\, E(i)$ of the injective left $C$-comodule $E(i)$. Then (a) follows by applying Lemma 2.5 (a) to $M= E(i)$.

(c) By (a), the $j$th column of $\Car_C$ is the $j$th row of $\Car_{C^{op}} $. Hence, (c) follows from (b) applied to the coalgebra $C^{op}$.

(d) By (b), the matrix $\Car_C$ is row finite if and only if $\dim_KE(i)$ is finite, for any $i\in I_C$, or equivalently, if and only if the injective envelope of any simple left $C$-comodule is of  finite $K$-dimension. But this property is equivalent to the right semiperfectness of $C$, see \cite{Lin}. 

(e) According to (c), the matrix $\Car_C$ is column  finite if and only if   the injective envelope of any simple right $C$-comodule is of  finite $K$-dimension. Since   this property is equivalent to the left semiperfectness of $C$  \cite{Lin}, the statement  (e) follows and the proof is complete. 
\end{proof}
\medskip

In \cite{Sim4}, a class of computable coalgebras  $C$, called left Euler coalgebras,  is  defined in such a way that the left Cartan matrix $\Car_C \in \mathbb{M}_{I_C}(\mathbb{Z})$ of such a coalgebra $C$ has a left inverse ${\Large{\mbox{$\mathfrak{c}$}}}_C^{-}$    in the non-associative  matrix algebra  $ \mathbb{M}_{I_C}(\mathbb{Z})$. Unfortunately, usually a left inverse ${\Large{\mbox{$\mathfrak{c}$}}}_C^{-}$ is not row-finite or column-finite. Below, we introduce a class of Euler coalgebras $C$ such that the left inverse ${\Large{\mbox{$\mathfrak{c}$}}}_C^{-}$ of $\Car_C$ is a   row-finite and a  column-finite matrix. 

We would like to remark here  that the multiplication in the   matrix algebra $ \mathbb{M}_{I_C}(\mathbb{Z})$ is not associative.  The matrices in  $ \mathbb{M}_{I_C}(\mathbb{Z})$ may have unequal left and right inverses and that one-sided inverse of a matrix may not be unique. Moreover, the left inverse may exist, without a right inverse existing. Also being invertible  as a $\mathbb{Z}$-linear map is not equivalent to being invertible as a matrix, see \cite{WZ}.
\medskip

We introduce a class of left Euler coalgebras as follows.\medskip

  {\sc Definition 2.7.} A $K$-coalgebra $C$ is defined to be a \textbf{left} (resp. right)  \textbf{sharp  Euler coalgebra} if  $C$ has  the following two properties.
  
  (a)  $C$ is computable, that is, $\dim_K\Hom_C(E', E'') $ is finite, for 
every pair of indecomposable injective left $C$-comodules $E'$ and $E''$. 
  
  (b) Every simple left (resp. right) $C$-comodule $S$ admits a  finite and socle-finite injective resolution
 \vspace{-1.2ex}
$$
0\,\,\mapr{}{}\,\, S  \,\,\mapr{}{}\,\,  E _0
\,\,\mapr{h_1}{}\,\, E _1
\,\,\mapr{h _2}{}\,\,\ldots 
\,\,\mapr{h _n}{}\,\,E _n  \,\,\mapr{}{}\,\,
0, 
  \vspace{-1.2ex}  
$$
 that is,  the injective comodules  $E_0, \ldots, E_n$ are  socle-finite.
 
 A $K$-coalgebra  $C$ is defined to be a \textbf{sharp Euler coalgebra} if it  is both left and right sharp Euler coalgebra and the following condition is satisfied
  
  (c) $\dim_K\Ext^m_C(S , S') = \dim_K\Ext^m_{C^{op}}(\widehat{S'}, \widehat{S}) $, for  all $m\geq 0$ and all simple left $C$-comodules $S$ and $S'$,   where  $\widehat{S} = DS $  and $\widehat{S'} = DS' $ are the dual simple right $C$-comodules.\bigskip
 
 Obviously, any sharp Euler coalgebra is an Euler coalgebra in the sense of \cite{Sim4}.
  Now we show that   (one-sided) semiperfect coalgebras $C$ with $\textrm{gl.dim}\, C<\infty$ are  sharp Euler coalgebras.\medskip
 
 {\sc Lemma 2.9.} \textit{Assume that $C$ is a pointed   left or right semiperfect coalgebra,  with a fixed decomposition $\soc\,{}_CC = \bigoplus _{j\in I_C} S(j)$.}
 
 (a)   $\Ext^m_C(S(a), S(b)) \cong   \Ext^m_{C^{op}}(\widehat{S}(b), \widehat{S}(a)) $, \textit{for  all $m\geq 0$ and  any pair of  simple left $C$-comodules  $S(a)$ and $S(b)$, with 
   $a,b\in I_C$, where $\widehat{S}(b) = DS(b)$ and $\widehat{S}(a)= DS(a)$  are the dual simple right $C$-comodules corresponding to   $a,b\in I_C$.}
 
 (b)  \textit{If the   global dimension} $\textrm{gl.dim}\, C$  \textit{of $C$ is  finite, then $C$ is a  sharp Euler coalgebra}.\medskip

 \begin{proof}  We prove the lemma in case $C$ is a pointed   left semiperfect coalgebra. The proof in case $C$ is   right semiperfect follows in a  similar way. 
 
 (a) Let   $S(a)$ and $S(b)$ be simple left $C$-comodules. Since $C$ is   left semiperfect,  there is a minimal projective resolution $\textbf{P}_*(a)$ of  $S(a)$ in $C\comod$. By Lemma 1.2, there is a duality $D: C\comod\,\,\mapr{}{}\,\, C^{op}\comod$ that   carries  $\textbf{P}_*(a)$ to a minimal injective resolution $D\textbf{P}_*(a)$ of $\widehat{S}(a)$ in the category $C^{op}\comod$ and induces an isomorphism of chain complexes $\Hom_C(\textbf{P}_*(a), S(b)) \cong \Hom_{C^{op}}(\widehat{S}(b), D\textbf{P}_*(a))$. Hence, we get the  induced isomorphism  of the cohomology $K$-spaces
 $$
 \Ext^m_C(S(a), S(b))= H^m[\Hom_C(\textbf{P}_*(a), S(b))] \cong H^m[\Hom_{C^{op}}(\widehat{S}(b), D \textbf{P}_*(a))] =  \Ext^m_{C^{op}}(\widehat{S}(b), \widehat{S}(a)) ,
 $$
 and (a) follows.
 
 (b) Assume that $\textrm{gl.dim}\, C$  of $C$ is  finite. Since $C$ is left semiperfect,  the indecomposable injectives in $C^{op}\Comod$ are finite-dimensional and therefore any simple right $C$-comodule has a finite injective resolution  in $C^{op}\comod$. Hence, $C$ is  a  right sharp Euler coalgebra. To prove that $C$ is   a left   sharp Euler coalgebra, assume that $S$ is a left simple $C$-comodule and let 
 \vspace{-1.2ex}
$$
0\,\,\mapr{}{}\,\, S  \,\,\mapr{}{}\,\,  E _0
\,\,\mapr{h_1}{}\,\, E _1
\,\,\mapr{h _2}{}\,\,\ldots 
\,\,\mapr{h _n}{}\,\,E _n  \,\,\mapr{}{}\,\,
0, 
  \vspace{-1.0ex}  
$$
be a minimal  injective  resolution of $S$ in $C\Comod$. We show that 
   the injective comodules  $E^{(j)}_0, \ldots, E^{(j)}_n$ are socle-finite. Assume that $S= S(b)$. Since a minimal projective resolution $\textbf{P}_*(b)$ of  $S=S(b)$ lies in $C\comod$ and is of  finite length $ \leq   \textrm{gl.dim}\, C$, it follows that,  for $m\geq 0$,  $\dim_K\Ext^m_C(S(a), S(b)) $ is finite, for all $a\in I_C$, and $\dim_K\Ext^m_C(S(a), S(b)) =0$, for all but a finite number of simple comodules $S(a)$.  Since  $C$ is pointed,  $\dim_K\Ext^m_C(S(a), S(b)) $ is the Bass number $\mu_m(S(a), S(b)) $ of the pair $(S(a), S(b)) $, that is, $\mu_m(S(a), S(b)) $ is the multiplicity the indecomposable injective comodule $E(a)$ appears in $E_m$, as a direct summand, see \cite[(4.23)]{Sim5}. It follows that, for each $m\geq 0$, the number $\mu_m(S(a), S(b)) $ is finite,  and $\mu_m(S(a), S(b)) $ is non-zero,  for at most finitely many $m$ and a finite number of indices $a\in I_C$.  Consequently, 
   the injective comodules  $E^{(j)}_0, \ldots, E^{(j)}_n$ are socle-finite, and the proof is complete.
 \end{proof}\medskip
 
  Next we give a description of sharp Euler  path  coalgebras $C= K^\SQR\!  Q$, with  $Q$ a quiver.\medskip
 
 {\sc Lemma 2.10.} \textit{Assume that $Q$ is a connected quiver and  $K^\SQR\!  Q$ is path  $K$-coalgebra  of a quiver  $Q$. The following three conditions are equivalent.}
 
 (a)  $K^\SQR\!  Q$  \textit{is  a sharp  Euler coalgebra}. 
 
 (b) $K^\SQR\!  Q$  \textit{is left and right Euler coalgebra.}
 
 (c) \textit{The quiver $Q$ is locally finite, that is, every vertex of $Q$ has at most finitely many neighbours in} $Q$. 
  \medskip

 \begin{proof}  The equivalence of (b) and (c) follows from \cite[Theorem 5.1(a)]{Sim4} and the implication (a)$\Rightarrow$(b) is obvious.  Since the coalgebras $C=K^\SQR\!  Q$ and  $C^{op}= (K^\SQR\! Q) ^{op} \cong  K^\SQR\!  Q^{op}$  are hereditary then to prove the inverse implication  (b)$\Rightarrow$(a),  it is enough to show that there is a $K$-linear isomorphism $\Ext^1_C(S(a), S(b)) \cong   \Ext^1_{C^{op}}(\widehat{S}(b), \widehat{S}(a)) $, for    any pair of  simple left $C$-comodules  $S(a)$ and $S(b)$, with 
   $a,b\in Q_0$, where $\widehat{S}(b) = DS(b)$ and $\widehat{S}(a)= DS(a)$  are the dual simple right $C$-comodules corresponding to   $a,b\in Q_0$.       Since the elements of $\Ext^1_C(S(a), S(b)) $ can be interpreted as equivalence classes of one-fold extensions $0\to S(b) \to N \to S(a) \to 0$ in $C\comod$ and  the duality $D:C\comod \to C^{op}\comod$ carries $0\to S(b) \to N \to S(a) \to 0$ to the exact sequence $0\to \widehat{S}(a) \to DN \to \widehat{S}(b) \to 0$,  it defines a $K$-linear  isomorphism $\Ext^1_C(S(a), S(b)) \cong   \Ext^1_{C^{op}}(\widehat{S}(b), \widehat{S}(a)) $.  This finishes the proof.
    \end{proof}\medskip
 
 Now we give   examples of   non-semiperfect sharp Euler coalgebras of infinite global dimension and of arbitrary large finite global dimension.\medskip
 
 {\sc Example 2.11.} Let $I$ be the infinite poset of the form
 \vspace{-0.4cm}
$$
 \vcenter{\vbox{ \mfpic[10][10]{-12}{6}{-2}{2}
  \tlabel[cc](-20,0){$\cdots$} 
    \tlabel[cc](20,0){$\cdots$} 
 \tlabel[cc](-12,0){$\bullet$} 
  \tlabel[cc](-6,0){$\bullet$} 
\tlabel[cc](0,0){$\bullet$} 
\tlabel[cc](-15,0){$\mathcal{I}_{-2}$} 
\tlabel[cc](-9,0){$\mathcal{I}_{-1}$} 
\tlabel[cc](-3,0){$\mathcal{I}_1$} 
\tlabel[cc](3,0){$\mathcal{I}_2$} 
\tlabel[cc](9,0){$\mathcal{I}_3$} 
\tlabel[cc](15,0){$\mathcal{I}_4$} 
\tlabel[cc](-12.0,1.1){$\scriptstyle{-1}$}
\tlabel[cc](-6.0,1.1){$\scriptstyle{0}$}
\tlabel[cc](0.0,1.1){$\scriptstyle{1}$}
\tlabel[cc](6,0){$\bullet$}
\tlabel[cc](6.0,1.1){$\scriptstyle{2}$}
\tlabel[cc](12,0){$\bullet$}
\tlabel[cc](12.0,1.1){$\scriptstyle{3}$}
\tlabel[cc](3,2){}
\tlabel[cc](3,0.5){}
 \tlabel[cc](3,-0.5){}
\tlabel[cc](3,-2){} 
{\headshape{2}{2}{false}
\arrow\arc{(-17.7,-0.4),(-12.3,-0.4),90} 
\arrow\arc{(-17.7,0.4),(-12.3,0.4),-90} 
\arrow\arc{(-11.7,-0.4),(-6.3,-0.4),90} 
\arrow\arc{(-11.7,0.4),(-6.3,0.4),-90} 
\arrow\arc{(-5.7,-0.4),(-0.3,-0.4),90} 
\arrow\arc{(-5.7,0.4),(-0.3,0.4),-90} 
\arrow\arc{(0.3,-0.4),(5.7,-0.4),90} 
\arrow\arc{(0.3,0.4),(5.7,0.4),-90}
\arrow\arc{(6.3,-0.4),(11.7,-0.4),90} 
\arrow\arc{(6.3,0.4),(11.7,0.4),-90}
\arrow\arc{(12.3,-0.4),(17.7,-0.4),90} 
\arrow\arc{(12.3,0.4),(17.7,0.4),-90}}
\endmfpic
}} \vspace{-0.4cm}
$$
directed from the left to the right,  where $\mathcal{I}_{m}= \mathcal{G}_{m}$ is the garland of length $|m| +1$ (see \cite{Si92}) \vspace{-0.2cm}
$$
\scriptsize{\begin{matrix}{\mathcal{G}_{m}: \quad \bullet\hspace{-0.5em}
\begin{array}[c]{c}
\nearrow\vspace{-0.2ex}\\
\searrow
\end{array}}
\hspace{-1.3em}
\begin{array}[c]{lcr}
\bullet\rightarrow\bullet- & \cdots & \rightarrow\bullet\rightarrow
\bullet\vspace{-0.3ex}\\
\hspace*{0.8em}{\nearrow\hskip-1.05em\searrow} & &
{\nearrow\hskip
-1.05em\searrow}\hspace*{1em}
{\nearrow\hskip-1.05em\searrow}\hspace*{0pt}\vspace{-0.3ex}\quad\\
\bullet\rightarrow\bullet- & \cdots & \rightarrow\bullet\rightarrow\bullet
\end{array}
\hspace{-1.2em}
{
\begin{array}[c]{c}
\searrow\vspace{-0.2ex}\\
\nearrow
\end{array}
\hspace{-0.4em}\ast }\,\,\,\qquad\mbox{$(2|m|+2$ vertices, $|m|\geq 1)$}.
\end{matrix}}\vspace{-0.2cm}
$$
 Obviously,  $I$ is an intervally finite poset. By the results given in \cite{Sim3} and \cite{Sim6}, the incidence coalgebra  $C = K^\SQR\! I$  of the poset $I$  has the following properties (see also \cite[Examples 4.25 and 4.26]{Sim4}):
 
 $1^\circ$  $C$ is  a sharp Euler coalgebra  and the   global dimension $\textrm{gl.dim}\, C$ of $C$ is infinite. 

 $2^\circ$  If $S(a)$ is the simple left $C$-comodule  corresponding to the vertex $a$ then the   injective  dimension $\textrm{inj.dim}\, S(a)$ of $S(a)$ equals $|m|+1$, for any $m\neq 0$. 
 
$3^\circ$   $\dim_K\Ext^m_C(S(a), S(b)) = \dim_K\Ext^m_{C^{op}}(\widehat{S}(b), \widehat{S}(a)) $, for  all $m\geq 0$ and all $a,b\in I$, where $\widehat{S}(b) = DS(b)$ and $\widehat{S}(a)= DS(a)$  are the simple right $C $-comodules corresponding to the vertices $a$ and $b$ in $I$.
 
$4^\circ$  $C$ is  both left and right locally artinian, locally cocoherent,   and the category  $C\Comod_{fc}$ coincides with the full subcategory of $C\Comod$ consisting of artinian objects.    

 $5^\circ$ The coalgebra  $C$ is neither left semiperfect nor right semiperfect. 

 $6^\circ$ The Cartan $\mathbb{Z} \times \mathbb{Z}$ square matrix $\Car_C \in  \mathbb{M}_{\mathbb{Z}}(\mathbb{Z})$   of $C$ is lower triangular and  has no finite rows and  no finite  columns.  
 
  $7^\circ$  $\Car_C$ has a unique    left inverse  $\Car^-_C \in  \mathbb{M}_{\mathbb{Z}}(\mathbb{Z})$, which is also a unique right inverse of $\Car_C$.  The matrix $\Car^-_C$ is row-finite and column-finite.

$8^\circ$ Let $m_0\geq 1$ be a fixed   integer and let $\mathcal{G}_{m_0}$ be the garland of length $|m_0| +1$. If we take $\mathcal{I}_{m}=  \mathcal{G}_{m_0}$, for each $m\in  \mathbb{Z}$,   in the construction of $I  $ then  $C= K^\SQR\! I$ is  a sharp Euler coalgebra,    $\textrm{gl.dim}\, C = m_0+1$  is finite,  $C$ is neither left semiperfect nor right semiperfect, and satisfies the conditions $4^\circ$, $6^\circ$, and $7^\circ$.
\medskip

  Now, given a pointed sharp Euler coalgebra $C$,  we construct a left inverse and a right inverse of the Cartan matrix $\Car_C \in \mathbb{M}_{I_C}(\mathbb{Z})$. We follow the proof of Theorem 4.18 in \cite{Sim4}, and the notation introduced there. Given a left (resp. right) sharp Euler coalgebra $C$, we fix  a    finite  minimal injective resolution  \vspace{-0.2cm}
  $$
0\,\,\mapr{}{}\,\, S(j) \,\,\mapr{h^{(j)}_0}{}\,\,  E^{(j)}_0
\,\,\mapr{h^{(j)}_1}{}\,\, E^{(j)}_1
\,\,\mapr{h^{(j)}_2}{}\,\,\ldots 
\,\,\mapr{h^{(j)}_n}{}\,\,E^{(j)}_n  \,\,\mapr{}{}\,\,
0, 
  \vspace{-1.2ex} \leqno(2.12)
$$
of the simple left $C$-comodule $S(j)$  in $C\Comod_{fc}$, with $E^{(j)}_0 = E(j)$,  and a    finite  minimal injective resolution  \vspace{-0.2cm}
  $$
0\,\,\mapr{}{}\,\, \widehat{S}(j) \,\,\mapr{\widehat{h}^{(j)}_0}{}\,\,  \widehat{E}^{(j)}_0
\,\,\mapr{\widehat{h}^{(j)}_1}{}\,\, \widehat{E}^{(j)}_1
\,\,\mapr{\widehat{h}^{(j)}_2}{}\,\,\ldots 
\,\,\mapr{\widehat{h}^{(j)}_{\widehat{n}}}{}\,\,\widehat{E}^{(j)}_{\widehat{n} } \,\,\mapr{}{}\,\,
0, 
  \vspace{-1.2ex} \leqno(2.13)
$$
of the simple right $C$-comodule $\widehat{S}(j)= DS(j)$ in $C^{op}\Comod_{fc}$,  with $\widehat{E}^{(j)}_0 = \widehat{E}(j)$, respectively. We fix    finite direct sum decompositions  \vspace{-0.7ex}
$$
E^{(j)}_m= \bigoplus _{p\in I_C} E(p)^{d ^{(j)}_{mp}}= \bigoplus _{p\in I^{(j)}_m} E(p)^{d ^{(j)}_{mp}}  ,\quad \widehat{E}^{(j)}_m= \bigoplus _{p\in I_C} \widehat{E}(p)^{\widehat{d} ^{(j)}_{mp}}= \bigoplus _{p\in \widehat{I}^{(j)}_m} \widehat{E}(p)^{\widehat{d} ^{(j)}_{mp}}  \vspace{-2.0ex}\leqno(2.14)
$$
of $E^{(j)}_m$ and  $\widehat{E}^{(j)}_m$,  for $m\geq 0$, where  $I^{(j)}_m$ and $\widehat{I}_m^{(j)}$ are  a finite subsets of $I_C$,  ${d} ^{(j)}_{mp}$ and $ \widehat{d} ^{(j)}_{mp}$ are  a  positive integers, for each $p\in I^{(j)}_m$ and each $p\in \widehat{I}^{(j)}_m$, respectively, and we set $d ^{(j)}_{mp}=0$, for any $p\in I_C\setminus I^{(j)}_m$, and $\widehat{d} ^{(j)}_{mp}=0$, for any $p\in {I}_C\setminus \widehat{I}^{(j)}_m$.   \medskip

 {\sc  Theorem 2.15.} \textit{Let $C$ be a pointed   computable  $K$-coalgebra,  with a fixed decomposition $\soc\,{}_CC = \bigoplus _{j\in I_C} S(j)$,  and let $\Car_C =[ \mathbf{c}_{ij}]_{i,j\in I_C} \in \mathbb{M}_{I_C}(\mathbb{Z})$ be  the left Cartan matrix} \textrm{(2.4)}  \textit{of $C$}.
   
  \textrm{(a)}   \textit{If   $C$ is a left sharp Euler coalgebra then the matrix $\Car_C^\leftarrow =[ \mathbf{c}^-_{ij}]_{i,j\in I_C} \in \mathbb{M}_{I_C}(\mathbb{Z})$, with $ \mathbf{c}^-_{jp} = 
\sum\limits_{m=0}^\infty (-1)^m  d^{(j)}_{mp} \in \mathbb{Z}
 $, is row-finite and is a left inverse of $\Car_C$ in $ \mathbb{M}_{I_C}(\mathbb{Z})$, where  $d^{(j)}_{mp} $ is the integer defined by the decomposition} (2.14)  \textit{of the $m$th term  $E^{(j)}_m$ of the minimal injective resolution} (2.12)  \textit{of the simple left $C$-comodule $S(j)$. Moreover, for each $j\in I_C$, we have}  $ \bdim\,  \widehat{E}(j) \cdot (\Car_C^\leftarrow)^{tr} = \bdim\ S(j)= e_j$, where $\bdim\, \widehat{E}(j)$ is the $j$th column of $\Car_C$.    
  
   \textrm{(b)}   \textit{If   $C$ is a right sharp Euler coalgebra then the matrix $\Car_C^\rightarrow =[ \widehat{\mathbf{c}}^-_{ij}]_{i,j\in I_C} \in \mathbb{M}_{I_C}(\mathbb{Z})$, with $ \widehat{\mathbf{c}}^-_{jp} = 
\sum\limits_{m=0}^\infty (-1)^m  \widehat{d}^{(p)}_{mj} \in \mathbb{Z}
 $, is column-finite and is a right inverse of $\Car_C$ in $ \mathbb{M}_{I_C}(\mathbb{Z})$, where  $\widehat{d}^{(p)}_{mj} $ is the integer defined by the decomposition} (2.14)  \textit{of the $m$th term  $\widehat{E}^{(p)}_j$ of the minimal injective resolution} (2.13) \textit{of the simple right $C$-comodule $\widehat{S}(p) = DS(p)$. Moreover, for each $j\in I_C$, we have}  $ \bdim\, E(j) \cdot  \Car_C^\rightarrow  = \bdim\ S(j)= e_j$. 
   
    \textrm{(c)}   \textit{If $C$ is a   sharp Euler coalgebra then  the matrix \vspace{-0.2cm}
    $$
    \Car_C^{-1}:=  \Car_C^\leftarrow = \Car_C^\rightarrow = [ \mathbf{c}^-_{ij}]_{i,j\in I_C} \in \mathbb{M}_{I_C}(\mathbb{Z}) ,\vspace{-0.2cm}\leqno(2.16)
    $$
     with $ \mathbf{c}^-_{ij}= \widehat{\mathbf{c}}^-_{ij}=  
\sum\limits_{m=0}^\infty (-1)^m   {d}^{(j)}_{mi} =  
\sum\limits_{m=0}^\infty (-1)^m  \widehat{d}^{(i)}_{mj} $,  is both row-finite and column-finite,  and  $\Car_C^-$ is  a left inverse of $\Car_C$ and  a right inverse of $\Car_C$.  } \medskip
 
\begin{proof} 
(a)   Assume   that $C$ is a left sharp Euler coalgebra. Then  the minimal injective  resolution (2.13) of $S(j)$ is finite and  the injective comodules  $E^{(j)}_0, \ldots, E^{(j)}_n$ are socle-finite. Hence the  sum
  $
 \mathbf{c}^-_{jp}=
\sum\limits_{m=0}^\infty (-1)^m  d^{(j)}_{mp}  
 $ 
 is an integer,  the matrix  $\Car_C^\leftarrow =[ \mathbf{c}^-_{ij}]_{i,j\in I_C} $  is well defined, and each of its row is finite, because the set $I^{(j)}_0 \cup I^{(j)}_1 \cup\ldots\cup  I^{(j)}_n \subseteq I_C$ is finite and  
 
\hspace{1cm} $\mathbf{c}^-_{jp}=
\sum\limits_{m=0}^\infty (-1)^m  d^{(j)}_{mp}= 
\sum\limits_{m=0}^n (-1)^m  d^{(j)}_{mp} = 0 $, for all $ p\not\in I^{(j)}_0 \cup I^{(j)}_1 \cup\ldots\cup  I^{(j)}_n $.  

To prove the   equality $\Car_C^\leftarrow\cdot  \Car_C = \textbf{E}$ (the identity matrix), we note that, by the additivity of the function $\bdim$,      the exact sequence (2.12) together with the decomposition (2.14) yields\vspace{-0.2cm}
$$
e_j= \bdim\, S(j)  = \sum\limits_{m=0}^{\infty} (-1)^m\bdim\, E^{(j)}_m = \sum\limits_{m=0}^\infty (-1)^m \sum _{p\in I_C} d ^{(j)}_{mp}\cdot \bdim\, E(p) = \sum _{p\in I_C}   \mathbf{c}^-_{jp}\cdot \bdim\, E(p) ,
\vspace{-1.5ex} 
$$
Hence the equality $\Car_C^\leftarrow\cdot  \Car_C =  \textbf{E}$ follows, because the $p$th row of the matrix $\Car_C$ is the dimension vector $\textbf{e}(p)= \bdim\, E(p)$ of $E(p)$, see (2.4).

By applying the matrix transpose $(-)^{tr}:\mathbb{M}_{I_C}(\mathbb{Z})\to \mathbb{M}_{I_C}(\mathbb{Z})$ and the equality $\Car_{C^{op}} =  \Car^{tr}_C$, we get $ \textbf{E}=  \textbf{E}^{tr} =  \Car_C^{tr} \cdot (\Car_C^\leftarrow) ^{tr}  = \Car_{C^{op}}  \cdot (\Car_C^\leftarrow) ^{tr} $ and, in view of Lemma 2.6, the equality $ \bdim\,  \widehat{E}(j) \cdot (\Car_C^\leftarrow)^{tr} = \bdim\ S(j)= e_j$ follows.

(b)  Assume   that $C$ is a right sharp Euler coalgebra.  Then  $C^{op}$ is a left sharp Euler coalgebra  and, by   (a)  with $C$  and $C^{op}$ interchanged,   the matrix ${\Car}_{C^{op}}^\leftarrow =[ \widehat{\mathbf{c}}^\leftarrow_{ij}]_{i,j\in I_C} \in \mathbb{M}_{I_C}(\mathbb{Z})$, with $  {\mathbf{c}}^\leftarrow_{jp} = 
\sum\limits_{m=0}^\infty (-1)^m  \widehat{d}^{(j)}_{mp} \in \mathbb{Z}
 $, is row-finite and is a left inverse of $\Car_{C^{op}} =  \Car^{tr}_C$  in $ \mathbb{M}_{I_C}(\mathbb{Z})$, where  $ \widehat{d}^{(j)}_{mp} $ is the integer defined by the decomposition (2.14)   of the $m$th term  $ \widehat{E}^{(j)}_m$ of the minimal injective resolution of the simple right $C$-comodule $ \widehat{S}(j) = DS(j)$.   It follows that $  {\mathbf{c}}^\leftarrow_{jp} =  \widehat {\mathbf{c}}^-_{pj}$, for all $j,p\in I_C$, and consequently, we get $({\Car}_{C^{op}}^\leftarrow)^{tr} = \Car_C^\rightarrow $

By  Lemma 2.6, we get  $\Car_{C^{op}} =  \Car^{tr}_C$ and the $p$th row $ \bdim\,  \widehat{E}(j)$ of  $\Car_{C^{op}}$  is the $p$th column of $\Car_C$. Since the equality  ${\Car}_{C^{op}}^\leftarrow \cdot \Car_{C^{op}}  =  \textbf{E}$ holds,   the matrix transpose yields $ \textbf{E} =  \textbf{E}^{op} =\Car^{tr}_{C^{op}} \cdot ({\Car}_{C^{op}}^\leftarrow)^{tr}  =   \Car_C\cdot \Car_C^\rightarrow$, that is, $\Car_C^\rightarrow$ is a right inverse of $\Car_C$. Hence (b) follows. 

(c)  Assume that $C$ is a sharp Euler coalgebra, that is, $C$ is left and right sharp  and the equality  $\dim_K\Ext^m_C(S(a), S(b)) = \dim_K\Ext^m_{C^{op}}(\widehat{S}(b), \widehat{S}(a)) $ holds, for all $a,b\in I_C$. We show  that $\Car_C^\leftarrow = \Car_C^\rightarrow $. Since $C$ is pointed, we have  \vspace{-0.2cm}
$$
  {d}^{(j)}_{mi} =\dim_K\Ext^m_C(S(i), S(j))\quad \mbox{and}\quad   \widehat{d}^{(i)}_{mj}= \dim_K\Ext^m_{C^{op}}(\widehat{S}(i), \widehat{S}(j)), \vspace{-0.2cm}\leqno(2.17)
  $$
see \cite[(4.23)]{Sim5}, and therefore   ${d}^{(j)}_{mi}=  \widehat{d}^{(i)}_{mj}$. It follows that, given $i,j\in I_C$, we have  $  \mathbf{c}^-_{ji} = 
\sum\limits_{m=0}^\infty (-1)^m  d^{(j)}_{mi} =  \sum\limits_{m=0}^\infty (-1)^m   \widehat{d}^{(i)}_{mj} =   \widehat{ \mathbf{c}}^-_{ji} $. This shows that $ \Car_C^\leftarrow=  \Car_C^\rightarrow$ and,  according to  (a) and (b), the matrix $\Car_C^{-1}:= \Car_C^\leftarrow=  \Car_C^\rightarrow$ (2.16) is   row-finite and column-finite, and is both left and right inverse of $\Car_C$.  This finishes the proof of the theorem.
\end{proof} 
\smallskip

{\sc Corollary 2.18}. \textit{Assume that $C$ is a pointed    sharp Euler coalgebra as in Theorem} 2.15,  \textit{with the Cartan matrix $\Car_C$ and its   inverse $\Car_C^{-1}$} (2.16). 

(a)  \textit{The matrix  $\Car^{-1}_C$ is row-finite and column-finite, and, given $a\in I_C$, we have:} 

 $\qquad \bdim\, E(a) \cdot \Car_C^{-1} =  \bdim\, S(a)$ \textit{and} $ \bdim\, S(a) \cdot \Car_C   =  \bdim\, E(a)$,

 $\qquad \bdim\,  \widehat{E}(a) \cdot (\Car_C^{-1})^{tr} =  \bdim\, S(a)$   \textit{and} $  \Car_C \cdot ( \bdim\, S(a))^{tr} =  (\bdim\, \widehat{E}(a))^{tr}$. 

(b) \textit{The subsets} $\{\bdim\, E(a)\}_{a\in I_C}$, $\{\bdim\, \widehat{E}(a)\}_{a\in I_C}$  \textit{of the group $\mathbb{Z}^{I_C}$ are $\mathbb{Z}$-linearly independent}.

(c) \textit{For each $j\in I_C$, the vector} $e_j = \bdim\, S(j) $ \textit{belongs to the subgroup generated by the  set } $\{\bdim\, E(a)\}_{a\in I_C}$, \textit{and to the subgroup generated by the  set} $\{\bdim\, \widehat{E}(a)\}_{a\in I_C}$.  \medskip

\begin{proof}  The equalities in (a) follow from Theorem 2.15, and (b) is a consequence of (a), because the vectors $\bdim\, S(a)= e_a\in  \mathbb{Z}^{I_C}$, with $a\in I_C$,     are $\mathbb{Z}$-linearly independent.

(c) We recall that the $a$th row of $\Car_C$ is the vector $\bdim\, E(a)$. Since the matrix $\Car_C^{-1}$ is row-finite,  the equality  $\Car^{-1}_C\cdot \Car_C= \textbf{E}$ yields $  e_j= \bdim\, S(j)   = \sum _{a\in I_C}   \mathbf{c}^-_{ja}\cdot \bdim\, E(a)$  \hbox{and the first part} of (c) follows. The second one follows in a similar way from the equality  $\Car_C\cdot \Car^{-1}_C= \textbf{E}$.  \end{proof}  \medskip

{\sc Corollary 2.19}. \textit{If  $C$ is   pointed    and  left semiperfect} (\textit{resp.  right semiperfect})  of finite global dimension \textit{then   the Cartan matrix $\Car_C$ of $C$ is column-finite} (\textit{resp. row-finite})  \textit{and the matrix}  $\Car_C^{-1}$  (2.16) \textit{is a two-sided inverse of $\Car_C$.  Moreover, $\Car_C^{-1}$} \textit{is column-finite and  row-finite, and the equalities of Corollary} 2.18(a) \textit{hold}. \medskip

\begin{proof}  By Lemma 2.6 (d) and  (e), the Cartan matrix $\Car_C$ of $C$ is column-finite (resp. row-finite), if  $C$ is  left semiperfect (resp.  right semiperfect).   Since, according to   Lemma 2.9,   $C$ is a  sharp Euler coalgebra,  the corollary follows from Corollary 2.18. \end{proof} 

  \vspace{-0.3cm}
 
 \section{Coxeter transformation for a sharp  Euler coalgebra } 
  \vspace{-0.3cm}
  
$\quad$  We study in this section  the properties of the Coxeter  
transformations  defined in \cite[Definition 4.27]{Sim4} for pointed   Euler coalgebras. Here, we also follow   \cite[Definition III.3.14]{ASS}. We modify \cite[Definition 4.27]{Sim4}  as follows.\medskip


{\sc Definition 3.1.}  Assume that   $C$ is a pointed    sharp 
 Euler  
$K$-coalgebra with fixed decomposition $\soc\,{}_CC = \bigoplus _{j\in I_C} S(j)$. Let  $  \Car_C \in \mathbb{M}_{I_C}(\mathbb{Z})$  be  the Cartan matrix  of $C$ and let $\Car_C^{-1}$   be   the  two-sided inverse  (2.16) of $\Car_C$.

 (a) The {\bf Coxeter matrix} of $C$ is the  $I_C\times I_C$ square matrix 
 $
 {\Phi}_C = - \Car_C^{-tr}\cdot \Car_C $, where we set $\Car_C^{-tr} =( \Car_C^{-1})^{tr} = ( \Car_C^{tr})^{-1}$.

(b) The  {\bf   Coxeter transformations} of $C$  are  the group homomorhisms\vspace{-0.3cm}
 $$
   \mathbb{Z}^{ I_C }_\blacktriangleright \,\,\RightLeftarr{\mathbf{\Phi}_C}{\mathbf{\Phi}^-_C}\,\, \mathbb{Z}^{I_C}_\blacktriangleleft \vspace{-0.3cm}\leqno(3.2)
   $$ 
 defined by the formulas
  $\mathbf{\Phi}_C(x) =  - (x\cdot \Car_C^{-tr})\cdot \Car_C $,
for $x\in  \mathbb{Z}^{I_C}_\blacktriangleright$, and  $\mathbf{\Phi}^-_C(y) =  - (y\cdot \Car_C^{-1})\cdot \Car^{tr}_C $,
for $y\in  \mathbb{Z}^{I_C}_\blacktriangleleft$, where $  \mathbb{Z}^{I_C}_\blacktriangleright   \subseteq \mathbb{Z}^{I_C}$ is  the subgroup  of  $\mathbb{Z}^{I_C}$ generated by the subset      $\{\widehat{\textbf{e}}(a) =\bdim\, \widehat{E}(a)\}_{a\in I_C}$  and $    \mathbb{Z}^{I_C}_\blacktriangleleft \subseteq \mathbb{Z}^{I_C}$ is  the subgroup  of  $\mathbb{Z}^{I_C}$ generated by the subset   $\{\textbf{e}(a) =\bdim\, E(a)\}_{a\in I_C}$.\bigskip

 By Corollary 2.18, the sets $\{\bdim\, E(a)\}_{a\in I_C}$ and  $\{\bdim\, \widehat{E}(a)\}_{a\in I_C}$ are $\mathbb{Z}$-linearly independent in $\mathbb{Z}^{I_C}$ and therefore they form   $\mathbb{Z} $-bases of  $  \mathbb{Z}^{(I_C)}_\blacktriangleleft$ and $  \mathbb{Z}^{(I_C)}_\blacktriangleright $, respectively.  Note also that $M\mapsto \bdim\, M$ defines  the group  isomorphism  of the Grothendieck group $K^\blacktriangleright_0(C) = K_0(C^{op}\inj)$ and $   \mathbb{Z}^{(I_C)}_\blacktriangleright  $, and  the group  isomorphism  of the Grothendieck group $K^\blacktriangleleft_0(C) = K_0(C\inj)$ and $  \mathbb{Z}^{I_C}_\blacktriangleleft$. Note that, by Corollary 2.18,  
 
 $\mathbf{\Phi}_C(\bdim\,  \nabla_CE(a)) = \mathbf{\Phi}_C(\widehat{\textbf{e}}(a)) =  - (\widehat{\textbf{e}}(a)\cdot \Car_C^{-tr})\cdot \Car_C = -e_a \cdot \Car_C = - {\textbf{e}}(a) = -   \bdim\,  E(a) $, 
 
  $\mathbf{\Phi}^-_C( \bdim\,  E(a))= \mathbf{\Phi}^-_C(\textbf{e}(a)) =  -e_a \cdot \Car_C^{-1}= -\widehat{\textbf{e}}(a)  - \bdim\,  \nabla_CE(a)$. 
  
  \noindent It follows that the transformations (3.2) are well-defined and mutually   inverse.  
 
 The following  theorem  is  the main result of this section (compare with  \cite[Corollary IV.2.9]{ASS}).  \smallskip
 
 {\sc Theorem  3.4.}  \textit{Assume that  $C$ is a pointed    sharp 
 Euler  
$K$-coalgebra with fixed decomposition $\soc\,{}_CC = \bigoplus _{j\in I_C} S(j)$. Let
$\mathbf{\Phi}_C$ and $\mathbf{\Phi}^-_C$ be the Coxeter transformations $(3.2)$ of $C$}.

(a)  Let $M$ be an
indecomposable   left  $C$-comodule in $C \textrm{-Comod}^\bullet_{fc}$ such that  $\operatorname*{inj.dim}\, M = 1$ and
$\operatorname*{Hom}_{C}(C,M)=0$. If  \vspace{-0.2cm}
$$  
0  \, \mapr{}{} \,\   M   \,\mapr{}{}  \,M' \,\mapr{}{} \,\tau^{-}_CM  \,
\mapr{}{} \, 0  \vspace{-0.2cm}
$$ 
 \textit{is the unique almost split sequence} (1.21)
\textit{in} $C\Comod_{fc}$, \textit{with an indecomposable comodule $\tau^{-}_CM$ lying in} $C \textrm{-comod}_{f\mathcal{P}}$ then  \vspace{-0.2cm}
$$
 \bdim\,\tau^{-}_CM=\mathbf{\Phi}^-_C(\bdim\, M).\vspace{-0.2cm}
$$

(b)    \textit{Assume that $N$ is  an
indecomposable non-projective   left  $C$-comodule in} $C \textrm{-comod}_{f\mathcal{P}}\subseteq C \textrm{-Comod} _{fc}$  \textit{such that    $\operatorname*{inj.dim}\, DN  =1$ and $\operatorname*{Hom} _{C}(C,DN)=0$.  If }\vspace{-0.2cm}
$$  
0  \, \mapr{}{} \,\  \tau_C N   \,\mapr{}{}  \,N' \,\mapr{}{} \,N  \,
\mapr{}{} \, 0  \vspace{-0.2cm}    
$$ 
\textit{is the unique almost split sequence} (1.22) \textit{in} $C\Comod_{fc}$, \textit{with an indecomposable comodule $\tau _CN$ lying in} $C \textrm{-Comod}^\bullet_{fc}$, then  
 $$
  \bdim\,\tau_CN = \mathbf{\Phi}_C(\bdim \, N).
$$
 
\begin{proof} (a)    Assume that $M$ is  an
indecomposable   left  $C$-comodule in $C \textrm{-Comod}^\bullet_{fc}$ such that  $\operatorname*{inj.dim}\, M = 1$. Then $M$ admits a minimal 
injective copresentation \vspace{-0.2ex}
$$
0\, \mapr{}{}\,M\, \mapr{}{}\,E_{0}\, \mapr{g}{}\, E_{1} \, \mapr{}{}\,0\vspace{-0.2ex}
$$
 in $C\Comod_{fc}$  where $E_{0}$ and $E_{1}$ are
socle-finite injective comodules.  It follows that $\bdim\,M = \bdim \, E_1 - \bdim \, E_0$. Since $\mathbf{\Phi}^-_C(\bdim\, E(a)) =  -\bdim\, E(a) \cdot \Car_C^{-1}= -\bdim\,\widehat{E}(a) $, for every $a\in I_C$ and the comodules  $E_{0}$ and $E_{1}$ are finite direct sums of the comodules ${E}(a) $, with  $a\in I_C$,   we get  $\mathbf{\Phi}^-_C(\bdim\, E_0) =  -\bdim\, \nabla_C({E}_0)$, $\mathbf{\Phi}^-_C(\bdim\, E_1) =  -\bdim\, \nabla_C({E}_1)$ and, by applying $\mathbf{\Phi}^-_C$,  the equality $\bdim\,M = \bdim \, E_1 - \bdim \, E_0$ yields  \vspace{-0.2ex}
$$
\mathbf{\Phi}^-_C(\bdim\, M) = \mathbf{\Phi}^-_C(\bdim \, E_1) - \mathbf{\Phi}^-_C(\bdim \, E_0)  =
 \bdim\, \nabla_C({E}_0)  -\bdim\, \nabla_C({E}_1) .  \vspace{-0.2ex}
$$
On the other hand,  the  exact sequence (1.10) in  $C^{op} \textrm{-Comod}$, induced by the injective copresentation of $M$, has the form  \vspace{-0.2ex}
$$
  0  \, \mapr{}{}\,   \textrm{Tr}_C(M)   \, \mapr{}{}\,  \nabla_C(E_1)\, \mapr{ \nabla_C(g)}{}\,  \nabla_C(E_0)  \, \mapr{}{}\, 0,\vspace{-0.2cm} 
  $$
  because the assumption $\operatorname*{Hom}_{C}(C,M)=0$ yields  $ \nabla_C(M)= \operatorname*{Hom}_{C}(C,M)^\circ = 0$, see Theorem 1.8 (a).  Since $\dim_K \textrm{Tr}_C(M) $ is finite, we have    $ \bdim\, D \textrm{Tr}_C(M) =   \bdim\,   \textrm{Tr}_C(M) $ and the exact sequence yields   \vspace{-0.2ex}
  $$
   \bdim\, \tau^-_C\, M =  \bdim\, D \textrm{Tr}_C(M) =   \bdim\,   \textrm{Tr}_C(M)  =   \bdim\, \nabla_C({E}_0)  -\bdim\, \nabla_C({E}_1)  = \mathbf{\Phi}^-_C(\bdim\, M)   \vspace{-0.2ex}
   $$
   and (a) follows.

 (b)  By Proposition 1.16 (b), there is a duality $D: C \textrm{-comod}_{f\mathcal{P}}  \, \, \mapr{\simeq}{}\, \,C^{op} \textrm{-comod}_{fc} $ that carries the indecomposable left $C$-comodule in $  C \textrm{-comod}_{f\mathcal{P}} $ to the indecomposable right $C$-comodule in $C^{op} \textrm{-comod}_{fc} $.  Since we assume $\operatorname*{inj.dim}\, DN  =1$,  the comodule $DN$ is not injective and there is a minimal socle-finite injective copresentation  \vspace{-0.2ex}
 $$ 
  0  \, \mapr{}{}\,   DN   \, \mapr{}{}\, E'_0 \, \mapr{g'}{}\,  E'_1   \, \mapr{}{}\, 0 \vspace{-0.2ex}
  $$ 
of $DN$ in $C^{op} \textrm{-Comod}_{fc} $.   By an obvious $\nabla_{C^{op}}$ version of Theorem 1.8, there is a short exact sequence
 $$
  0  \, \mapr{}{}\,   \textrm{Tr}_{C^{op}}(DN)   \, \mapr{}{}\, \nabla_{C^{op}}(E'_1)\, \mapr{\nabla_{C^{op}}(g')}{}\, \nabla_{C^{op}}(E'_0)  \, \mapr{}{}\, \nabla_{C^{op}}(DN) \, \mapr{}{}\,0 
  $$
 in $C\Comod_{fc} $. The assumption   ${\Hom} _{C}(C,DN)=0$ yields   $\nabla_{C^{op}}(DN) = \Hom_C(C, DN)^\circ = 0$. Then, by applying the arguments used in the proof of (a),  we get   \vspace{-0.2ex}
  $$
  \begin{array}{rcl}
  \bdim\,   \textrm{Tr}_{C^{op}}(DN) &= &\bdim \,  \nabla_{C^{op}}(E'_1) -  \bdim \,  \nabla_{C^{op}}(E'_0)\\
   &= & \mathbf{\Phi} _C(\bdim \,   E'_0) -  \mathbf{\Phi} _C(\bdim \,    E'_1)\\
   &= &  \mathbf{\Phi} _C(\bdim\, DN) =  \mathbf{\Phi} _C(\bdim\, N) , 
   \end{array}  \vspace{-0.2ex}
$$
because $\dim_KN$ is finite. This finishes the proof.
\end{proof} 
  \medskip

   {\sc Remark  3.5.}  If $C$ is a sharp Euler  coalgebra such that $\textrm{gl. dim}\,  C=1$ and $M$ (resp. $N$)  is an indecomposable non-injective comodule (resp. non-projective comodule), we have   $\textrm{inj. dim}\,  M=1$ and $\operatorname*{Hom}_ C(C, M)=0$ (resp.  $\textrm{inj. dim}\,  DN=1$ and $\operatorname*{Hom}_ C(C, DN)=0$), and Theorem 3.4 applies to $M$ (resp. to  $N$). \medskip


 \section{Illustrative examples}
 
 In this section we illustrate previous results by concrete examples.\medskip

 {\sc Example  4.1.} Let   $Q$ be the  infinite locally Dynkin quiver  \vspace{-1.6ex}
$$
Q :   \begin{array}{l}
 \vspace{-1.3ex} \\
{\bullet}
 \, \mapr{}{}  \,
{\bullet}
 \, \mapr{}{}  \, {\bullet}  \, \mapr{}{} \,  {\bullet}
\, \mapr{}{} \, {\bullet}
\, \mapr{}{}  \,  {\bullet }
\, \mapr{}{}   
\ldots \vspace{-1.5ex} \\
{\scriptstyle 0}
 \, \qquad   
{\scriptstyle 1}
 \,  \qquad  {\scriptstyle 2}  \,  \qquad   {\scriptstyle 3}
  \qquad  {\scriptstyle 4}\,
   \qquad    {\scriptstyle 5 }
\,  
\end{array}\vspace{-0.2ex}
 $$
 of type $\mathbb{A}_\infty$ and let $C=   K^\SQR\!  Q$ be the path $K$-coalgebra of $Q$, see \cite{Ch}, \cite{Sim1}, \cite{Wo}. Then $C$ has the upper triangular matrix form 
 $$
C =\left[\scriptsize{\begin{matrix}
  &K &K & K & K& K& K&K&K&\ldots\cr
  &0&K & K & K& K& K&K&K&\ldots\cr
  &0&0 & K & K& K& K&K&K&\ldots\cr
  &0&0 & 0 &  K& K& K&K&K&\ldots\cr
  &0&0 & 0& 0& K& K&K&K&\ldots\cr
  &0&0& 0& 0& 0& K&K&K&\ldots\cr
  &0&0&0& 0& 0& 0& K&K&\ldots\cr
  &0&0&0& 0& 0& 0& 0&K&\ldots\cr
  & \vdots & \vdots & \vdots& \vdots& \vdots&\vdots &\vdots
&\vdots   &\ddots &\cr
\end{matrix}}
\right]
\vspace{-0.2ex}
$$ 
and consists of the upper triangular $\mathbb{N}\times \mathbb{N}$ square matrices with coefficients in $K$ with at most  finitely many 
non-zero entries. Then $\soc{}_CC =   \bigoplus _{j\in I_C} S(j)$, where $I_C= \mathbb{N} = \{0,1,2,\ldots\}$ and $S(n) = Ke_n$ is the simple subcoalgebra spanned by the matrix $e_n\in C$ with $1$ in the $ n\times n$ entry, and zeros elsewere. Note that $e_n$ is a group-like element of $C$.

 The Cartan  matrix $\Car_C\in \mathbb{M} _{\mathbb{N}}(\mathbb{Z})$    of
$C$ and its   inverse $\Car_C^{-1}$ have the lower triangular forms \vspace{-0.2ex}
$$
\Car_C =\left[\scriptsize{\begin{matrix}
  &1 &0 & 0 & 0& 0& 0&0&0&\ldots\cr
  &1&1& 0 & 0& 0& 0&0&0&\ldots\cr
  &1&1 &1 & 0& 0& 0&0&0&\ldots\cr
  &1&1 &1 &1& 0& 0&0&0&\ldots\cr
  &1&1 &1& 1&1& 0&0&0&\ldots\cr
  &1&1&1&1&1&1&0&0&\ldots\cr
  &1&1&1&1&1&1&1&0&\ldots\cr
  &1&1&1&1&1&1&1&1&\ldots\cr
  & \vdots & \vdots & \vdots& \vdots& \vdots&\vdots &\vdots
&\vdots   &\ddots &\cr
\end{matrix}}
\right]\qquad \Car^{-1}_C =\left[\scriptsize{\begin{matrix}
  &1 &0 & 0 & 0& 0& 0&0&0&\ldots\cr
  &\!\!\!-1&1& 0 & 0& 0& 0&0&0&\ldots\cr
  &0&\!\!\!-1 &1 & 0& 0& 0&0&0&\ldots\cr
  &0&0 &\!\!\!-1 &1& 0& 0&0&0&\ldots\cr
  &0&0&0&\!\!\! -1&1& 0&0&0&\ldots\cr
  &0&0&0&0&\!\!\!-1&1&0&0&\ldots\cr
  &0&0&0&0&0&\!\!\!-1&1&0&\ldots\cr
  &0&0&0&0&0&0&\!\!\!-1&1&\ldots\cr
  & \vdots & \vdots & \vdots& \vdots& \vdots&\vdots &\vdots
&\vdots   &\ddots &\cr
\end{matrix}}
\right]
\vspace{-0.2ex}
$$ 
Hence the Coxeter matrices $ {\Phi} _C= - \Car^{-tr}_C\cdot \Car_C$ and ${\Phi}^{-1} _C= - \Car^{-1}_C\cdot \Car^{tr}_C$ are of the forms\vspace{-0.2ex}
$$
 {\Phi} _C=\left[\scriptsize{\begin{matrix}
  &0 &1 & 0 & 0& 0& 0&0&0&\ldots\cr
  &0&0& 1 & 0& 0& 0&0&0&\ldots\cr
  &0&0 &0 & 1& 0& 0&0&0&\ldots\cr
  &0&0&0&0& 1& 0&0&0&\ldots\cr
  &0&0&0&0&0& 1&0&0&\ldots\cr
  &0&0&0&0&0&0&1&0&\ldots\cr
  &0&0&0&0&0&0&0&1&\ldots\cr
  &0&0&0&0&0&0&0&0&\ldots\cr
  & \vdots & \vdots & \vdots& \vdots& \vdots&\vdots &\vdots
&\vdots   &\ddots &\cr
\end{matrix}}
\right]\qquad  {\Phi}^{-1}_C = \left[\scriptsize{\begin{matrix}
  &\!\!\!-1 &\!\!\!-1 & \!\!\!-1 & \!\!\!-1& \!\!\!-1& \!\!\!-1&\!\!\!-1&\!\!\!-1&\ldots\cr
  &1&0& 0 & 0& 0& 0&0&0&\ldots\cr
  &0&1 &0 & 0& 0& 0&0&0&\ldots\cr
  &0&0 & 1 &0& 0& 0&0&0&\ldots\cr
  &0&0&0& 1&0& 0&0&0&\ldots\cr
  &0&0&0&0& 1&0&0&0&\ldots\cr
  &0&0&0&0&0& 1&0&0&\ldots\cr
  &0&0&0&0&0&0& 1&0&\ldots\cr
  & \vdots & \vdots & \vdots& \vdots& \vdots&\vdots &\vdots
&\vdots   &\ddots &\cr
\end{matrix}}
\right]
\vspace{-0.2ex}
$$ 
The coalgebra $C$ is pointed, representation-directed in the sense of \cite{Sim4},  right semiperfect and hereditary, that is, $\textrm{gl. dim}\, C= 1$. Hence $C$ is a sharp  Euler  
$K$-coalgebra.   Every left $C$-comodule is a direct sum of finite-dimensional ones \cite{Sim1} and therefore every indecomposable left $C$-comodule is finite-dimensional.
The left $C$-comodules in $C\comod$ can be identified with the  finite-dimensional $K$-linear representations of the quiver $Q$. Under the identification $C\comod= \rep_K( Q)  $, the Auslander-Reiten quiver of  $C\comod$ has the form\vspace{-0.2cm}

\begin{center}
\underline{{\sc Figure 1.}  {\sc The Auslander-Reiten quiver of 
the category $ C\comod \cong  \rep_K( Q) $}}  \vspace{-0.2cm}
\end{center} 
\def\II{{\mathbb{I}}}
$$
\begin{array}{cccccccccccccccccc}
\ldots& _4{\II}_4 &&\hidewidth\hbox{ - - - - - - - - }\hidewidth&& _3{\II}_3 &
&\hidewidth\hbox{ - - - - - - - - }\hidewidth&& _2{\II}_2
&&\hidewidth\hbox{ -
- - - - - - - }\hidewidth&& _1{\II}_1  &&\hidewidth\hbox{ - - - - - - - -
}\hidewidth & &_0{\II}_0  \vspace{1.5ex}  \\
  &  &\searrow&   &\nearrow &    &\searrow &  &\nearrow &    &\searrow &&\nearrow
&    &\searrow && \nearrow&
 \vspace{1.5ex}  \\ 
 \ldots &\hidewidth\hbox{ - - - - - - - - }\hidewidth& & _3{\II}_4  &&
\hidewidth\hbox{ - - - - - - - - }\hidewidth & & _2{\II}_3  && \hidewidth\hbox{
- - - - - - - - }\hidewidth & & _1{\II}_2  &&\hidewidth\hbox{ - - - - - - - -
}\hidewidth&  &  _0{\II}_1  && \vspace{1.5ex}  \\
 & &\nearrow&    & \searrow &  &\nearrow &    &\searrow &&\nearrow
 &    &\searrow&& \nearrow&& 
&     \vspace{1.5ex} \\
 \ldots  &   _3{\II}_5 & & \hidewidth\hbox{ - - - - - - - - }\hidewidth &&
_2{\II}_4  && \hidewidth\hbox{ - - - - - - - - }\hidewidth && _1{\II}_3  &&
\hidewidth\hbox{ - - - - - - - - }\hidewidth  & &_0{\II}_2 & & && \vspace{1.5ex}  \\
       & &\searrow &&\nearrow
 &    &\searrow&& \nearrow&&\searrow
&    & \nearrow&&   &   &  &     \vspace{1.5ex}   \\
 \ldots   & \hidewidth\hbox{ - - - - - - - - }\hidewidth &&  _2{\II}_5 & &
\hidewidth\hbox{ - - - - - - - - }\hidewidth && _1{\II}_4  &&\hidewidth\hbox{ -
- - - - - - - }\hidewidth && _0{\II}_3  &&   & && &  
\\
        & &\nearrow
 &    &\searrow& & \nearrow& &\searrow
&   & \nearrow &&    &   &  &&&\\
    &\vdots& 
 &    & &\vdots&  & & 
&   \vdots & &&    &  &    &  &  & 
\end{array}\vspace{-0.2cm}
$$
  see \cite{NowSim} and \cite{Sim2}, where \vspace{-0.2cm}
  $$
_n{\II}_m  :  { 0}
 \, \mapr{}{}   \, \ldots  \, \mapr{}{} \,  { 0}
\, \mapr{}{} \, {  K_n}
\, \mapr{id}{}  \,  {  K_{n+1} }
\, \mapr{id}{}   
\ldots  \, \mapr{id}{}  \,  {\  K_m }\, \mapr{ }{}  \,    0 \,  \mapr{ }{} \,  0  \, \mapr{ }{}  \  \ldots \, \ldots \,  \vspace{-0.2cm}
 $$
 $K_n=K_{n+1}= \ldots = K_m = K$ and  $n\leq m$. Note that ${}_m{\II}_m= S(m)$ is simple and $_0{\II}_m  = E(m)$ is the injective envelope of $S(m)$, for each $m\geq 0$. Hence, the indecomposable injectives in the category $C\comod  $ form the right hand  section  \vspace{-0.2cm}
 $$
 \ldots \to {} _0{\II}_6\to  {}_0{\II}_5\to  {}_0{\II}_4\to  {}_0{\II}_3 \to  {}_0{\II}_2  \to  {}_0{\II}_1\to  {}_0{\II}_0 \vspace{-0.2cm}
 $$  
  of  Figure 1. Note also that  $C\comod $ 
 contains   no non-zero projective objects.   
Thus\vspace{-0.1cm}
$$
0 \, \mapr{}{}\,\operatorname{Tr}_C(_n{\II}_m) \, \mapr{}{}\, \nabla_CE(n-1) \, \mapr{}{}\,
 \nabla_CE(m-1)  \vspace{-0.2cm}
$$
is an injective  copresentation yielding  $\tau^-_C({}_n{\II}_m)= D\operatorname*{Tr}_C({}_n{\II}_m)\cong  {}_{n-1}{\II}_{m-1}$, for $n\geq 1$. \ The almost split sequences are\vspace{-0.2cm}
$$
0 \, \mapr{}{}\, {}_n{\II}_m \, \mapr{}{}\, {}_{n-1}{\II}_{m}\oplus {}_{n}{\II}_{m-1} \, \mapr{}{}\,
{}_{n-1}{\II}_{m-1} \, \mapr{}{}\,  0\vspace{-0.2cm}
$$
with irreducible morphisms ${}_{n}{\II}_{m} \, \mapr{}{}\, {}_{n}{\II}_{m+1}$ and ${}_{n-1}{\II}_{m} \, \mapr{}{}\,
{}_{n}{\II}_{m}$ being the obvious monomorphism into the first summand and epimorphism
onto the second summand. The map on the right is given by natural epimorphism
and monomorphism with alternate signs. This means that $\tau_C({}_{n-1}{\II}_{m-1}) \cong {}_n{\II}_m$ and  $\tau^-_C({}_n{\II}_m) \cong {}_{n-1}{\II}_{m-1}$, if $n\geq 1$.   Note also that 
$\bdim\, \tau_C({}_{n-1}{\II}_{m-1} )=  \mathbf{\Phi}_C(\bdim\, 
 {}_n{\II}_m)$ and  $\bdim \, \tau^-_C({}_n{\II}_m) =   \mathbf{\Phi}^{-1}_C(\bdim \ {}_{n-1}{\II}_{m-1})$,   if $n\geq 1$ (compare with Theorem 3.4).  \bigskip

  {\sc Example  4.2.} Let   $Q$ be the  infinite locally Dynkin quiver  \vspace{-1.6ex} 
  $$
Q :  \begin{array}{l}
 \vspace{-1.3ex} \\
\ldots \, \mapr{}{}  \,{\bullet}
 \,  \, \mapr{}{}  \,{\bullet}
 \,  \, \mapr{}{}  \,{\bullet}
 \, \mapr{}{}  \,
{\bullet}
 \, \mapr{}{}  \, {\bullet}  \, \mapr{}{} \,  {\bullet}
\, \mapr{}{} \, {\bullet}
\, \mapr{}{}  \,  {\bullet }
\, \mapr{}{}   
\ldots \vspace{-1.5ex} \\
 \, \ \ \qquad {\scriptstyle -2}
  \qquad  
\!{\scriptstyle -1}
 \, \qquad  {\scriptstyle 0}
 \, \qquad   
{\scriptstyle 1}
 \,  \qquad  {\scriptstyle 2}  \,  \qquad   {\scriptstyle 3}
  \qquad  {\scriptstyle 4}\,
   \qquad    {\scriptstyle 5 }
\,  
\end{array}
 $$
of type ${}_{ \infty} {\mathbb{A}}_{\infty}$  and let $C=   K^\SQR\!  Q$ be the path $K$-coalgebra of $Q$. Then $C$ has the upper triangular matrix form 
 $$
C =\left[\scriptsize{\begin{matrix}
\ddots  &\vdots &\vdots & \vdots & \vdots& \vdots& \vdots&\vdots&\vdots& \cr
\ldots  &K &K & K & K& K& K&K&K&\ldots\cr
\ldots   &0&K & K & K& K& K&K&K&\ldots\cr
\ldots  &0&0 & K & K& K& K&K&K&\ldots\cr
\ldots   &0&0 & 0 &  K& K& K&K&K&\ldots\cr
 \ldots  &0&0 & 0& 0& K& K&K&K&\ldots\cr
\ldots   &0&0& 0& 0& 0& K&K&K&\ldots\cr
 \ldots  &0&0&0& 0& 0& 0& K&K&\ldots\cr
 \ldots  &0&0&0& 0& 0& 0& 0&K&\ldots\cr
  & \vdots & \vdots & \vdots& \vdots& \vdots&\vdots &\vdots
&\vdots   &\ddots &\cr
\end{matrix}}
\right]
\vspace{-0.1ex}
$$ 
and consists of the upper triangular $\mathbb{Z}\times \mathbb{Z}$ square matrices with coefficients in $K$ with at most  finitely many 
non-zero entries. Then $\soc{}_CC =   \bigoplus _{j\in I_C} S(j)$, where $I_C= \mathbb{Z} = \{\ldots, -2,-1, 0,1,2,\ldots\}$ and $S(n) = Ke_n$ is the simple subcoalgebra spanned by the matrix $e_n\in C$ with $1$ in the $ n\times n$ entry, and zeros elsewere. Note that $e_n$ is a group-like element of $C$. The coalgebra $C$ is pointed, hereditary, left and right locally artinian and,  by Corollary 2.10,  $C$ is a sharp  Euler  
$K$-coalgebra.  Obviously, $C$ is neither right semiperfect nor left   semiperfect.

 The Cartan  matrix $\Car_C\in \mathbb{M} _{\mathbb{N}}(\mathbb{Z})$    of
$C$ and its   inverse $\Car_C^{-1}$ have the lower triangular forms 
$$
\Car_C =\left[\scriptsize{\begin{matrix}
\ddots  &\ddots &\vdots  & \vdots  & \vdots & \vdots & \vdots &\vdots &\vdots & \cr
\ldots  &1 &0 & 0 & 0& 0& 0&0&0&\ldots\cr
\ldots   &1&1& 0 & 0& 0& 0&0&0&\ldots\cr
 \ldots  &1&1 &1 & 0& 0& 0&0&0&\ldots\cr
  \ldots &1&1 &1 &1& 0& 0&0&0&\ldots\cr
 \ldots  &1&1 &1& 1&1& 0&0&0&\ldots\cr
 \ldots  &1&1&1&1&1&1&0&0&\ldots\cr
 \ldots  &1&1&1&1&1&1&1&0&\ldots\cr
 \ldots  &1&1&1&1&1&1&1&1&\ldots\cr
  & \vdots & \vdots & \vdots& \vdots& \vdots&\vdots &\vdots
&\vdots   &\ddots &\cr
\end{matrix}}
\right]\qquad \Car^{-1}_C =\left[\scriptsize{\begin{matrix}
\ddots  &\ddots &\ddots  &   & \vdots & \vdots & \vdots &\vdots &\vdots &\vspace{-2ex} \cr
  \ddots &1 &0 & 0 & 0& 0& 0&0&0&\ldots\vspace{-2ex}\cr
 \ddots  &\!\!\!-1&1& 0 & 0& 0& 0&0&0&\ldots\cr
  &0&\!\!\!-1 &1 & 0& 0& 0&0&0&\ldots\cr
\ldots   &0&0 &\!\!\!-1 &1& 0& 0&0&0&\ldots\cr
\ldots   &0&0&0&\!\!\! -1&1& 0&0&0&\ldots\cr
\ldots   &0&0&0&0&\!\!\!-1&1&0&0&\ldots\cr
  \ldots &0&0&0&0&0&\!\!\!-1&1&0& \vspace{-2ex}\cr
 \ldots  &0&0&0&0&0&0&\!\!\!-1&1& \ddots\cr
  & \vdots & \vdots & \vdots& \vdots& \vdots&\vdots &\ddots
&\ddots   &\ddots &\cr
\end{matrix}}
\right]
\vspace{-0.1ex}
$$ 
Hence the Coxeter matrices $  {\Phi} _C= - \Car^{-tr}_C\cdot \Car_C$ and $ {\Phi}^{-1} _C= - \Car^{-1}_C\cdot \Car^{tr}_C$ are of the forms
$$
 {\Phi} _C=\left[\scriptsize{\begin{matrix}
 \ddots  &\ddots &  & \vdots  & \vdots & \vdots & \vdots &\vdots &\vdots &\vspace{-2ex}  \cr
\ddots &0 &1 & 0 & 0& 0& 0&0&0&\ldots\cr
 \ldots &0&0& 1 & 0& 0& 0&0&0&\ldots\cr
 \ldots &0&0 &0 & 1& 0& 0&0&0&\ldots\cr
  \ldots&0&0&0&0& 1& 0&0&0&\ldots\cr
  \ldots&0&0&0&0&0& 1&0&0&\ldots\cr
 \ldots&0&0&0&0&0&0&1&0&\ldots\cr
 \ldots &0&0&0&0&0&0&0&1&  \vspace{-2ex}\cr
 \ldots &0&0&0&0&0&0&0&0&\ddots\cr
  & \vdots & \vdots & \vdots& \vdots& \vdots&&\ddots 
&\ddots   &\ddots &\cr
\end{matrix}}
\right]\qquad 
 {\Phi}^{-1}_C = \left[\scriptsize{\begin{matrix}
 \ddots  &\ddots &  & \vdots  & \vdots & \vdots & \vdots &\vdots &\vdots &\vspace{-2ex}  \cr
 \ddots &0 &0 & 0&0&0&0&0&0&\ldots\vspace{-2ex} \cr
\ddots  &1&0& 0 & 0& 0& 0&0&0&\ldots\cr
\ldots  &0&1 &0 & 0& 0& 0&0&0&\ldots\cr
\ldots &0&0 & 1 &0& 0& 0&0&0&\ldots\cr
\ldots  &0&0&0& 1&0& 0&0&0&\ldots\cr
\ldots  &0&0&0&0& 1&0&0&0&\ldots\cr
 \ldots &0&0&0&0&0& 1&0&0&  \vspace{-2ex}\cr
\ldots &0&0&0&0&0&0& 1&0&\ddots\cr
  & \vdots & \vdots & \vdots& \vdots& \vdots& & 
\ddots&\ddots   &\ddots &\cr
\end{matrix}}
\right]
\vspace{-0.1ex}
$$ 

 It is known that there is an equivalence of categories 
 $
K^\SQR\!     Q\comod \cong \rep_K(  Q) 
$ 
and we view it as an identification,  see \cite{CKQ}, \cite{Sim1} and  \cite[Proposition 3.3]{Sim7}. We recall from        \cite[Corollary  5.13]{Sim2} that any finite-dimensional $K$-linear representation  $N\in \rep_K( Q)$ of the infinite quiver $ Q$ restricts to  a representation of a finite convex linear quiver $Q^N = \mathbf{supp} (N)$ (the support of $N$) of the Dynkin type $\mathbb{A}_n$  and  is isomorphic to a    finite
interval representation of the form
 $$
_n{\II}_m  : \, \ldots \, \mapr{}{} \,  { 0}
 \, \mapr{}{}    \,  { 0}
 \, \mapr{}{}   \, \ldots  \, \mapr{}{} \,  { 0}
\, \mapr{}{} \, {  K_n}
\, \mapr{id}{}  \,  {  K_{n+1} }
\, \mapr{id}{}   
\ldots  \, \mapr{id}{}  \,  {\  K_m }\, \mapr{ }{}  0 \,  \mapr{ }{} \,  0  \, \mapr{ }{}  \  \ldots \,  \vspace{-0.2cm}
 $$
where  $- \infty<m\leq  t<\infty $ and $K_j=K$, for all $m\leq j \leq t$.  It is easy to see that  the indecomposable injective $K^\SQR\!    Q$-comodules  are infinite-dimensional. Hence  the category $C\comod$ contains no non-zero injective objects and no non-zero projective objects.  
 
  By Corollary 1.23 (see also   \cite{NowSim} and \cite[Section 6]{Sim2}), 
  every indecomposable object $N$ of $C\comod  $ has an almost split sequence in  $C\comod$ starting from $N$  and  has an almost split sequence in $C\comod $ terminating in $N$. Moreover the Auslander-Reiten translation quiver $\Gamma(C\comod )$ of the category $C\comod  $   has the form 
\begin{center}
\underline{{\sc Figure 2.} {\sc The Auslander-Reiten quiver of 
the category $ C\comod \cong  \rep_K( Q))$ } }
\end{center}
$$
 {\begin{array}{ccccccccccccccccccccc}
\ldots& \hidewidth_4{\II}_4\hidewidth &&\hidewidth\hbox{ - - - - - -  }\hidewidth&&\hidewidth _3{\II}_3 \hidewidth&
&\hidewidth\hbox{ - - - - - -   }\hidewidth&&\hidewidth _2{\II}_2\hidewidth
&&\hidewidth\hbox{ - - - - - -   }\hidewidth&& \hidewidth_1{\II}_1 \hidewidth &&\hidewidth\hbox{ - - - - - -  
}\hidewidth & &\hidewidth_0{\II}_0\hidewidth&&\hidewidth\hbox{ - - - - - -  
}\hidewidth & \medskip \\
  &  &\searrow&   &\nearrow &    &\searrow &  &\nearrow &    &\searrow &&\nearrow
&    &\searrow && \nearrow&  &\searrow && \nearrow
 \medskip \\ 
 \ldots &\hidewidth\hbox{   - - - }\hidewidth& &\hidewidth _3{\II}_4 \hidewidth &&
\hidewidth\hbox{ - - - - - -  }\hidewidth & & \hidewidth_2{\II}_3 \hidewidth && \hidewidth\hbox{
- - - - - -  }\hidewidth & &\hidewidth _1{\II}_2 \hidewidth &&\hidewidth\hbox{ - - - - - -  
}\hidewidth&  & \hidewidth _0{\II}_1 \hidewidth &&\hidewidth\hbox{ - - - - - -  
}\hidewidth&&\hidewidth  _{\scriptstyle{-}1}{\II}_0 \hidewidth &\,\,\hidewidth\hbox{    - - - }\hidewidth  \medskip \\
 & &\nearrow&    & \searrow &  &\nearrow &    &\searrow &&\nearrow
 &    &\searrow&& \nearrow&& 
\searrow&& \nearrow &   &\searrow  \medskip \\
 \ldots  &  \hidewidth _3{\II}_5 \hidewidth& & \hidewidth\hbox{ - - - - - -   }\hidewidth &&
\hidewidth_2{\II}_4 \hidewidth && \hidewidth\hbox{ - - - - - -  }\hidewidth && \hidewidth_1{\II}_3\hidewidth &&
\hidewidth\hbox{ - - - - - -  }\hidewidth  & &\hidewidth _0{\II}_2 \hidewidth& &\hidewidth\hbox{ - - - - - - 
}\hidewidth   &&\,\, \hidewidth _{\scriptstyle{-}1}{\II}_1\hidewidth&&\hidewidth\hbox{ - - - - - - 
}\hidewidth & \medskip \\
       & &\searrow &&\nearrow
 &    &\searrow&& \nearrow&&\searrow
&    & \nearrow&&\searrow&& \nearrow&&\searrow
&&\nearrow   \medskip  \\
 \ldots   & \hidewidth\hbox{      - - - }\hidewidth &&\hidewidth  _2{\II}_5 \hidewidth& &
\hidewidth\hbox{ - - - - - -   }\hidewidth &&\hidewidth _1{\II}_4\hidewidth  &&\hidewidth\hbox{ -
- - - - -   }\hidewidth &&\hidewidth _0{\II}_3 \hidewidth &&\hidewidth\hbox{ - - - - - -  
}\hidewidth   & &\hidewidth_{\scriptstyle{-}1}{\II}_2\hidewidth & &  \hidewidth\hbox{ - - - - - -  
}\hidewidth &&\hidewidth_{\scriptstyle{-}2}{\II}_1\hidewidth&\,\, \hidewidth\hbox{    - - - }\hidewidth
\\
        & &\nearrow
 &    &\searrow& & \nearrow& &\searrow
&   & \nearrow &&\searrow& & \nearrow& &\searrow&  & \nearrow& &\searrow\\
    &\vdots& 
 &    & &\vdots&  & & 
&   \vdots & &&    & \vdots &    &  &  & \vdots  &&&
\end{array}}\vspace{-0.2cm}
$$
Note that the Coxeter transformation  $\mathbf{\Phi} _C:\mathbb{Z}^{ I_C }_\blacktriangleright  \to  \mathbb{Z}^{ I_C }_\blacktriangleleft $ (3.2) extends to the isomorphism  $
 \mathbf{\Phi} _C: \mathbb{Z}^{\mathbb{Z}} \,\,\mapr{}{}\,\, \mathbb{Z}^{\mathbb{Z}}$   defined by the formula $ \mathbf{\Phi} _C(x) = x\cdot  \mathbf{\Phi} _C $. It  carries any vector  $x = (x_n) _{n\in \mathbb{Z} }\in  \mathbb{Z}^{ \mathbb{Z}}$ to the vector $ \mathbf{\Phi} _C(x) = \widehat{ x}=   (\widehat{ x}_n) _{n\in \mathbb{Z}}\in  \mathbb{Z}^{( \mathbb{Z})}$, with $\widehat{ x}_n= x_{n-1}$, for all $n\in \mathbb{Z}= I_C$. This means that $ \mathbf{\Phi} _C$ shifts any  vector $x\in \mathbb{Z}^{\mathbb{Z}}$ by one step to the right.  It follows that the inverse $
 \mathbf{\Phi} _C^{-1}: \mathbb{Z}^{\mathbb{Z}} \,\,\mapr{}{} \,\, \mathbb{Z}^{\mathbb{Z}}$ of $ \mathbf{\Phi} _C$ shifts any  vector $x\in \mathbb{Z}^{\mathbb{Z}}$ by one step to the left. 

Hence, by applying  the   the Auslander-Reiten quiver shown in Figure 2, we conclude that, given   an  indecomposable  $N$ in $C\comod$,    there
exist  almost split sequences
  \vspace{-1ex}
$$  
0  \, \mapr{}{} \,\  \tau_C N   \,\mapr{}{}  \,Y \,\mapr{}{} \,N  \,
\mapr{}{} \, 0\quad \mbox{and}\quad  0  \, \mapr{}{} \,\   N  \,\mapr{}{}  \,Z \,\mapr{}{} \,\tau_C^{-1}N \,
\mapr{}{} \, 0
 \vspace{-1ex}
$$ 
in  $C\comod  $  
and   the following equalities hold (compare with Theorem 3.4)
  \vspace{-1ex}
 $$
 \textbf{dim} (\tau_C N )=  \mathbf{\Phi} _C ( \textbf{dim}\, N)\quad \mbox{and}\quad  \textbf{dim} ( \tau_C^{-1} N) =  \mathbf{\Phi} _C^{-1} ( \textbf{dim}\,  N).
 \vspace{-1ex}
$$ 

Let us also look at the (abelian) category $C\Comod_{fc}$ of finitely copresented $C$-comodules.  It consists of artinian $C$-comodules and, by applying \cite[Proposition 2.13(a)]{Sim4}, one can show that  every indecomposable comodule $M$ of $C\Comod_{fc}$ is either injective, with $\dim_KM = \infty$, or $M$ is finite-dimensional  isomorphic to one of the comodules listed in Figure 2 and $\dim_K \mathrm{Tr}_C M $ is finite.  It follows that $C\Comod^\bullet_{fc}=  C\Comod_{fc}$ and $C\comod^\bullet_{fc}=  C\comod$. One can also show that the Grothendieck group $K_0(C\Comod_{fc})$ of $C\Comod_{fc}$  is isomorphic to the   
  Grothendieck group $K^\blacktriangleleft_0(C) = K_0(C\inj)\cong   \mathbb{Z}^{I_C}_\blacktriangleleft$
of the category $C\inj$.  Moreover,  the   Auslander-Reiten quiver $\Gamma(C\Comod_{fc})$ of the category $C\Comod_{fc}$ has two connected components: 

(a) the component shown in Figure 2  consisting of all indecomposable $C$-comodules of finite dimension, and 

(b) the following component consisting of all indecomposable injective $C$-comodules:
$$
\ldots \, \mapr{}{}\, E(-2)  \, \mapr{}{}\, E(-1)  \, \mapr{}{}\, E(0)  \, \mapr{}{}\, E(1)  \, \mapr{}{}\, E(2)  \, \mapr{}{}\, E(3)  \, \mapr{}{}\, \ldots \, .
$$
\smallskip
 
 {\sc Example  4.3.} Let   $Q= (Q_0, Q_1)$ be the  infinite locally Dynkin quiver of type $\mathbb{D}_\infty$ presented  in Example 1.25, with $Q_0=\{-1,0,1,2,3,\ldots \}$,  and let $C= K^\SQR\!  Q$ be the path $K$-coalgebra of $C$.  Then $C$ has a  $Q_0\times Q_0$ square matrix form shown in Example 1.25.  
The Cartan  matrix $\Car_C\in \mathbb{M} _{Q_0}(\mathbb{Z})$    of
$C$ and its   inverse $\Car_C^{-1}$ have the lower triangular forms \vspace{-0.2ex}
$$
\Car_C =\left[\scriptsize{\begin{matrix}
  &1&0 &1 & 0 & 0& 0& 0&0&0&  \ldots\cr
  &0&1&1& 0 & 0& 0& 0&0&0&  \ldots\cr
  &0&0&1&0&0 &0 & 0& 0& 0&  \ldots\cr
  &0&0&1&1 &0 &0& 0& 0&0& \ldots\cr
  &0&0&1&1 &1& 0&0& 0&0& \ldots\cr
  &0&0&1&1&1&1&0&0&0& \ldots\cr
  &0&0&1&1&1&1&1&0&0& \ldots\cr
  &0&0&1&1&1&1&1&1&0& \ldots\cr
  &0&0&1&1&1&1&1&1&1& \ldots\cr
  & \vdots& \vdots& \vdots & \vdots & \vdots& \vdots& \vdots&\vdots &\vdots
    &\ddots &\cr
\end{matrix}}
\right]\qquad \Car^{-1}_C =\left[\scriptsize{\begin{matrix}
  &1 &0 &\!\!\!-1 & 0& 0& 0&0&0&\ldots\cr
  &0&1&\!\!\!-1& 0 & 0& 0& 0&0& \ldots\cr
  &0&0 &1 & 0& 0& 0&0&0&\ldots\cr
  &0&0 &\!\!\!-1 &1& 0& 0&0&0&\ldots\cr
  &0&0&0&\!\!\! -1&1& 0&0&0&\ldots\cr
  &0&0&0&0&\!\!\!-1&1&0&0&\ldots\cr
  &0&0&0&0&0&\!\!\!-1&1&0&\ldots\cr
  &0&0&0&0&0&0&\!\!\!-1&1&\ldots\cr
  & \vdots & \vdots & \vdots& \vdots& \vdots&\vdots &\vdots
&\ddots   &\ddots &\cr
\end{matrix}}
\right]
\vspace{-0.2ex}
$$ 
Hence the Coxeter matrices $ {\Phi} _C= - \Car^{-tr}_C\cdot \Car_C$ and $  {\Phi}^{-1} _C= - \Car^{-1}_C\cdot \Car^{tr}_C$ are of the forms\vspace{-0.2ex}
$$
 {\Phi} _C=\left[\scriptsize{\begin{matrix}
  &\!\!\!-1 &0 & \!\!\!-1 & 0& 0& 0&0&0&\ldots\cr
  &0&\!\!\!-1& \!\!\!-1 & 0& 0& 0&0&0&\ldots\cr
  &1&1 &2 & 1& 0& 0&0&0&\ldots\cr
  &0&0&0&0& 1& 0&0&0&\ldots\cr
  &0&0&0&0&0& 1&0&0&\ldots\cr
  &0&0&0&0&0&0&1&0&\ldots\cr
  &0&0&0&0&0&0&0&1&\ldots\cr
  & \vdots & \vdots & \vdots& \vdots& \vdots&\vdots &\vdots
&\ddots   &\ddots &\cr
\end{matrix}}
\right]\qquad {\Phi}^{-1}_C = \left[\scriptsize{\begin{matrix}
  &0&1 &1 & 1& 1& 1&1&1&\ldots\cr
   &1&0 &1 & 1& 1& 1&1&1&\ldots\cr
  &\!\!\!-1&\!\!\!-1&\!\!\!-1& \!\!\!-1&\!\!\!-1& \!\!\!-1&\!\!\!-1&\!\!\!-1&\ldots\cr
  &1&1 & 1 &0& 0& 0&0&0&\ldots\cr
  &0&0&0&1& 0&0&0&0&\ldots\cr
    &0&0&0&0& 1&0&0&0&\ldots\cr
  &0&0&0&0&0& 1&0&0&\ldots\cr
  &0&0&0&0&0&0& 1&0&\ldots\cr
  & \vdots & \vdots & \vdots& \vdots& \vdots&\vdots &\vdots
&\ddots   &\ddots &\cr
\end{matrix}}
\right]
\vspace{-0.2ex}
$$ 
The coalgebra $C$ is pointed, right semiperfect and hereditary, that is, $\textrm{gl. dim}\, C= 1$. Hence $C$ is a sharp  Euler  
$K$-coalgebra.   Every left $C$-comodule is a direct sum of finite-dimensional ones \cite{Sim1} and therefore every indecomposable left $C$-comodule is finite-dimensional. If $N$ is a terminus of a mesh in the Auslander-Reiten quiver $\Gamma(C\comod)$ shown in Figure 0 then $N$ is the right hand term of an almost split sequence $  
0  \, \mapr{}{} \,\  \tau_C N   \,\mapr{}{}  \,N' \,\mapr{}{} \,N  \,
\mapr{}{} \, 0 $ in $C\comod$  
and    $
 \textbf{dim}(\tau_C N )=  \mathbf{\Phi} _C ( \textbf{dim}\,N)$, by Theorem    3.4.  It follows that the dimension vectors of the modules lying in each of two  infinite components of  $\Gamma(C\comod)$ shown in Figure 0 can be computed from the dimension vectors of the  modules lying on the sections $(*)$ and $(**)$ (presented in Example 1.25) by applying the iterations $  \mathbf{\Phi} _C^m $, with $m\geq 1$, of the Coxeter transformation $  \mathbf{\Phi} _C$.  Obviously, $C$ is representation-directed in the sense of \cite{Sim4} and every indecomposable $C$-comodule $N$ is uniquely determined by its dimension vector $ \textbf{dim}\,N$, see \cite{Sim1}-\cite{Sim4}.
 


\medskip

Addresses: 

 W. Chin,  Department of Mathematics, DePaul University,  Chicago, Illinois 60614,  U.S.A.
 
 e-mail: wchin@condor.depaul.edu 

 D. Simson,  Faculty of Mathematics and Computer Science, Nicolaus Copernicus University,  
 
 87-100 Toru\'n, ul. Chopina 12/18, Poland,  e-mail:  simson@mat.uni.torun.pl

\closegraphsfile

\small

\begin{thebibliography}{9999}                                                                                              

\bibitem{ASS}I. Assem, A. Skowro\'{n}ski and D. Simson,  \textit{Elements of the 
Representation Theory of Associative Algebras, Volume 1: Techniques of
Representation Theory}, Londom Math Soc. Student Texts
65, Cambridge University Press, Cambridge 2006.\vspace{-1ex}

\bibitem{Ausl}  M. Auslander, Coherent functors, \textit{in } Proc. Conf. on
Categorical Algebra, La Jolla, Springer-Verlag, 1966, pp. 189--231.\vspace{-1ex}


\bibitem{ARS}M. Auslander, I. Reiten, S. Smal\o ,  \textit{Representation
Theory of Artin algebras}. Cambridge Studies in Advanced Mathematics, 36.
Cambridge University Press, Cambridge 1997.\vspace{-1ex}



\bibitem{Ch}{}W. Chin, A brief introduction to coalgebra representation
theory, in: Proceedings from an International Conference Held at DePaul
University, J. Bergen, S. Catoiu, W. Chin, eds.   Lecture  Notes in
Pure and Appl. Math.,  Marcel-Dekker, 237(2004),  pp. 109--131.\vspace{-1ex}

 

\bibitem{CKQ}W. Chin, M. Kleiner and D. Quinn, Almost split sequences for
comodules, J. Algebra 249(2002), 1-19.\vspace{-1ex}


\bibitem{CKQ2}{}W. Chin, M. Kleiner and D. Quinn, Local theory of almost split
sequences for comodules, Ann. Univ. Ferrara - Sez. VII - Sc. Mat.Vol.  51(2005),
183-196.\vspace{-1ex}


\bibitem{CMo}W. Chin and S. Montgomery, Basic Coalgebras, \textit{in}
Modular Interfaces (Riverside, CA, 1995) AMS/IP Studies in Advanced Math.
vol.4, Providence RI, 4(1997), 41-47.\vspace{-1ex}


\bibitem{CG}J. Cuadra, J. G\'{o}mez-Torrecillas, Idempotents and
Morita-Takeuchi theory. Comm. Algebra 30 (2002),   2405--2426.\vspace{-1ex}


 \bibitem{DNR} S. D\u{a}sc\u{a}lescu, C. N\u{a}st\u{a}sescu  and
S. Raianu,   "Hopf Algebras. An Introduction",  {\em  Lecture Notes in
Pure and Applied Mathematics}, No. 235,  Marcel-Dekker, New-York,
2001.\vspace{-1ex}

\bibitem{Ga73}{}
P. Gabriel, Indecomposable representations II,    Symposia Mat.
Inst. Naz. Alta Mat.  11(1973), 81--104.\vspace{-1ex}


\bibitem{Gr}J.A. Green, Locally finite representations, J. Algebra 41(1976), 137-171.\vspace{-1ex}

 \bibitem{KR}    M. Kleiner and I. Reiten, Abelian categories, almost split
sequences and  comodules,  Trans. Amer. Math. Soc.  357(2005), 
 3201--3214.
 \vspace{-1.0ex}

\bibitem{Lin}B. I.-P. Lin, Semiperfect coalgebras, J. Algebra 49 (1977),  357-373.\vspace{-1ex}


\bibitem{Montg}  S. Montgomery,   \textit{Hopf Algebras and Their Actions on
Rings}, CMBS No. 82, AMS,  1993.\vspace{-1ex}


\bibitem{NTZ}C. Nastasescu, B. Torrecillas, Y. H. Zhang, Hereditary
coalgebras. Comm. Algebra 24(1996),   1521--1528.\vspace{-1ex}

\bibitem{NowSim} S. Nowak and  D. Simson, Locally Dynkin quivers and hereditary
coalgebras whose left comodules are direct sums of finite dimensional comodules,
Comm. Algebra 30(2002), 455-476.\vspace{-1ex}



\bibitem{Sim0}D. Simson, On the structure of pure semisimple Grothendieck
categories, Cahiers de Topologie et Geom. Diff. 33(1982), 397-406.\vspace{-1ex}


\bibitem{Si92}  D. Simson,    \textit{Linear Representations of Partially Ordered
Sets and Vector Space Categories},  Algebra, Logic and
Applications Vol. 4, Gordon \& Breach Science Publishers,
1992.\vspace{-1ex}


\bibitem{Sim1}{}
D. Simson, Coalgebras, comodules, pseudocompact algebras and  tame
comodule type,   Colloq. Math.   90(2001),
101--150.\vspace{-1.2ex}


\bibitem{Sim2}D. Simson, Coalgebras of tame comodule type, \textit{in:}
``Representations of Algebras'', Proceedings ICRA-9, (Eds. D. Happel and Y. B.
Zhang), Beijing Normal University Press, 2002, Vol. 2, pp. 450-486.\vspace{-1ex}


\bibitem{Sim3}D. Simson,  Path coalgebras of quivers with relations and a tame-wild
dichotomy problem for coalgebras, Lecture Notes in Pure Appl. Math.,  
Marcel Dekker 236(2004), pp. 465-492.\vspace{-1ex}

 \bibitem{Sim06}{}
D. Simson, Irreducible morphisms, the Gabriel quiver and colocalisations for 
coalgebras,  {\it Intern. J. Math. Math. Sci.,} 72(2006), 1--16.\vspace{-1ex} 

\bibitem{Sim4}
D. Simson, $\textrm{Hom}$-computable coalgebras,  a composition factors matrix and an Euler  bilinear
form of an Euler  coalgebra,   J. Algebra, 315(2007), 42--75.\vspace{-1ex}

 \bibitem{Sim5}{}
D. Simson, Localising embeddings of comodule    categories with applications to
tame and Euler coalgebras, J. Algebra, 312(2007), 455--494.
\vspace{-1ex}

 \bibitem{Sim7} D. Simson, Path coalgebras of  profinite bound quivers,   cotensor coalgebras of          bound species and   locally nilpotent
representations, {\em Colloq. Math.}  109(2007),  307--343.\vspace{-1ex}

 

 \bibitem{Sim08} D. Simson, Tame-wild dichotomy for  coalgebras, {\em J. London  Math. Soc }   78(2008),  783--797.\vspace{-1ex}
 
  \bibitem{Sim6}{}
D. Simson, Incidence   coalgebras  of intervally  finite posets, their integral quadratic forms and comodule categories, {\em Colloq. Math.}  115(2009),  259--295.
\vspace{-1ex}


\bibitem{Sut}R. Suter, Modules for $\emph{U} _{q}(sl_{2}),$ Comm. Math.
Phys., 162(1994), 359--393.\vspace{-1ex}


\bibitem{Tak}M. Takeuchi, Morita theorems for categories of comodules, J.
Fac. Sci. Univ. Tokyo 24(1977), 629--644.\vspace{-1ex}


\bibitem{Wo}D. Woodcock,   Some categorical remarks on the representation
theory of coalgebras, Comm. Algebra 25(1997), 2775--2794.\vspace{-1ex}

\bibitem{WZ} A. Wilansky and K. Zeller,  Inverses of matrices and matrix transformations,  Trans. Amer. Math. Soc. 6(1955), 414--420.

\end{thebibliography}
\end{document}